\documentclass[final,onefignum,onetabnum]{siamonline190516}

\pdfoutput=1

\usepackage[section]{placeins} 

\usepackage{algorithmic}
\usepackage[caption=false]{subfig}
\ifpdf
  \DeclareGraphicsExtensions{.eps,.pdf,.png,.jpg}
\else
  \DeclareGraphicsExtensions{.eps}
\fi

\usepackage{amssymb}
\usepackage{amsmath}
\usepackage{amsfonts}
\usepackage{upgreek}
\usepackage{siunitx}
\usepackage{listings}
\usepackage{multirow}
\usepackage{multicol}
\usepackage{bigdelim}
\usepackage{color}
\usepackage{stmaryrd}
\usepackage{soul}
\usepackage[symbol]{footmisc}
\usepackage{booktabs}
\usepackage[section]{placeins} 
\usepackage[titletoc]{appendix}
\usepackage{calc}
\usepackage[capitalize]{cleveref}
\usepackage[sort&compress,square,numbers]{natbib}
\usepackage{xspace}

\usepackage{graphicx}
\usepackage{graphics}
\usepackage{float}
\usepackage{wrapfig}
\usepackage[percent]{overpic}
\usepackage{varwidth}
\usepackage{epstopdf}
\usepackage{adjustbox}
\usepackage{tikz}
\usepackage{pgfplots}
\usetikzlibrary{arrows,matrix,positioning,fit}
\usetikzlibrary{shapes,positioning}
\usetikzlibrary{backgrounds}
\pgfplotsset{compat=1.14}
\usepgfplotslibrary{colorbrewer}
\usepgfplotslibrary{patchplots}
\usepgfplotslibrary[colorbrewer]
\usetikzlibrary{pgfplots.colorbrewer}
\usetikzlibrary[pgfplots.colorbrewer]
\usepgfplotslibrary{fillbetween}
\usepackage{tikz-layers}
\usepgfplotslibrary{units}
\usetikzlibrary{spy}
\usepackage{pgfplotstable}
\graphicspath{ {figs/} }
\usetikzlibrary{external}

\newcommand{\etal}{et~al.}
\newcommand\floor[1]{\lfloor#1\rfloor
\newcommand\ceil[1]{\lceil#1\rceil}}

\newcommand\hb[1]{\hat{\mathbf{#1}}}

\newcommand\tb[1]{\tilde{\mathbf{#1}}}
\newcommand\ob[1]{\overline{\mathbf{#1}}}

\newcommand\px[2]{\frac{\partial #1}{\partial {#2}}}

\newcommand\dx[2]{\frac{\mathrm{d} #1}{\mathrm{d} #2}}

\newcommand{\half}{\frac{1}{2}}

\newlength\myheight
\newlength\mydepth
\settototalheight\myheight{Xygp}
\newsiamremark{remark}{Remark}
\newsiamremark{hypothesis}{Hypothesis}
\crefname{hypothesis}{Hypothesis}{Hypotheses}
\newsiamthm{claim}{Claim}
\newsiamthm{prop}{Proposition}

\DeclareRobustCommand\onedot{\futurelet\@let@token\@onedot}
\def\@onedot{\ifx\@let@token.\else.\null\fi\xspace}

\def\ie{\emph{i.e}\onedot}

\def\etal{\emph{et al}\onedot}
\makeatother

\begin{document}

\title{An extended range of energy stable flux reconstruction methods on triangles\thanks{\corresponding{Will Trojak (\email{w.trojak@imperial.ac.uk})} \funding{None to declare.}}}
\headers{FR/SD on triangles}{W. Trojak and P. Vincent}

\author{Will Trojak\thanks{Department of Aeronautics, Imperial College London, South Kensington, London, SW7 2AZ}, \and Peter Vincent\footnotemark[2].}

\maketitle

\begin{abstract}
    We present an extended range of stable flux reconstruction (FR) methods on triangles through the development and application of the summation-by-parts framework in two-dimensions. This extended range of stable schemes is then shown to contain the single parameter schemes of \citet{Castonguay2011} on triangles, and our definition enables wider stability bounds to be developed for those single parameter families. Stable upwinded spectral difference (SD) schemes on triangular elements have previously been found using Fourier analysis. We used our extended range of FR schemes to investigate the linear stability of SD methods on triangles, and it was found that a only first order SD scheme could be recovered within this set of FR methods.
\end{abstract}



\section{Introduction}\label{sec:intro}
    The flux reconstruction (FR)~\citep{Huynh2007} method is a versatile numerical approach to approximating the solution of advection-diffusion equations and can be generalised to arbitrary order accuracy. A compelling advantage of the FR method is the dominance of locally structured computation in the algorithm, making it highly efficient on modern GPU hardware~\citep{Vincent2016}. Furthermore, the FR method can be thought of as a generalisation of the discontinuous Galerkin (DG)~\citep{Reed1973,Hesthaven2008} method over the set of different lifting functions. Within the FR literature, these lifting functions are realised through what are called a correction functions. In the seminal work of \citet{Huynh2007}, it was shown that changes to the correction function could result in significantly different numerical characteristics. 
    
    Families of correction functions can be formed by considering different norms in order to prove linear stability. The proofs of linear stability generally look to show that for a system such as:
    \begin{equation}\label{eq:lin_advec_1d}
        \px{u(x,t)}{t} + \px{f(u)}{x} = 0, \quad \mathrm{for} \quad u(x,t=0)=u_0(x), \quad x\in K\subset\mathbb{R}, \quad \mathrm{and} \quad f = au.
    \end{equation}
    The objective of stability proofs for FR was to find some correction function where the following is true:
    \begin{equation}
        \dx{}{t}\|u\|^2_A \leq 0,
    \end{equation}
    for some norm $A$, such as in the work of \citet{Vincent2010}. It was later shown by \citet{Allaneau2011} and \citet{Zwanenburg2016} that there is an equivalence between FR and linearly filtered DG. 
    The work of \citet{Vincent2015} defined a wider --- multi-parameter --- set of stable FR schemes which can be interpreted as a wider set of filters, with the filter implicitly defined by $A$. Concurrently, \citet{Ranocha2016} showed how FR could be cast into a summation-by-parts (SBP) formulation to achieve the same results.
    
    In the work of \citet{Castonguay2011}, and later that of \citet{Williams2013}, a stable set of FR schemes on triangles was defined that is analogous to those defined in 1D by \citet{Vincent2010}. For simulations of real flows, triangular elements are advantageous due to the ease of mesh generation for complex geometries, such as via Delaunay triangulation. Concurrently to these works on FR, \citet{Balan2011} defined a spectral difference~(SD) method on triangles. The FR and SD methods are closely related~\citep{Jameson2010}; however, previous works like that of \citet{Veilleux2022} have used Fourier analysis with upwinding to find stable SD schemes. In general, SD schemes on triangles are weakly unstable, and we would like to answer if stable SD methods on triangles can be found as a subset of linearly stable FR methods.
    
    The objective of this work is twofold, firstly we investigate the FR method on triangles further through the use of the SBP framework. We seek to extend the current definition of linearly stable FR methods on triangles from that of \citet{Castonguay2011} to a definition analogous to that of \citet{Vincent2015}. Furthermore, we will study the conditions needed for stability, symmetry, and conservation. Secondly, using this extended definition of stable FR on triangles we investigate the connection between FR and SD, and whether linearly stable SD methods can be found within this set of FR methods. With this in mind, this work is structured as follows. In \cref{sec:prelim} we introduce the FR method, previous correction functions, and the SBP framework. The main results of this work are presented in \cref{sec:linear}, where linear stability of the FR method on triangles is explored, and conditions are set out for conservation and symmetry. In \cref{sec:erfr_tri}, we define the new set of stable FR methods for several orders and prescribe there stability conditions. Then in \cref{sec:sd} we use this extended range of stable FR of schemes to investigate SD on triangles. In \cref{sec:numeical} we present some brief numerical results from using the newly found FR schemes and finally conclusions are drawn in \cref{sec:conclusions}.
\section{Preliminaries}\label{sec:prelim}
\subsection{Flux Reconstruction}
    The flux reconstruction (FR) method of \citet{Huynh2007} can be applied to advection and advection-diffusion systems as a method of spatial discretisation. In this work we solely focus on advection systems, and we now give a brief introduction to the FR method in one dimension. For a more detailed description of the method, readers can see references \citet{Huynh2007,Witherden2015} and references therein.
    
    Beginning with \cref{eq:lin_advec}, the first stage of FR is to subdivide the domain $K$ into $N$ compatible sub-domains, $\{K\}_{i\leq N}$, such that:
    \begin{equation}
        K  = \bigcup^N_{i=1} K_i \quad \mathrm{and} \quad K_i \cap K_j = \begin{cases}
            K_i &\mbox{if } i=j\\
            \emptyset &\mbox{otherwise}
        \end{cases}.
    \end{equation}
    A reference sub-domain $\hat{K}=[-1,1]$ is defined with the transformation $T_i:K_i\mapsto \hat{K}$, such that interpolation and differentiation operators can be more efficiently applied. Here we will only consider affine elements, \ie where $T_i$ is a linear functional. For a given polynomial order, $k$, $k+1$ solution points are placed in the sub-domains at the physical locations $\{x_{ij}\}_{0\leq j\leq k}$, and in the reference sub-domain at $\{\xi_{j}\}_{0\leq j\leq k}$. 
    
    Polynomials of the solution and flux functions from element $K_i$ can be fitted in $\hat{K}$ using a Lagrange finite-element as
    \begin{equation}
        \hat{u}_i(\xi) = \sum^k_{j=0}u_i(x_j)l_j(\xi), \quad \mathrm{and} \quad \hat{f}^\delta_i(\xi) = \sum^k_{j=0}f(u_i(x_j))l_j(\xi), \quad \mathrm{for}\quad l_n(\xi) = \prod^{k}_{\substack{j=0\\j\neq n}}\frac{\xi - \xi_j}{\xi_n - \xi_j}.
    \end{equation}
    To the flux we have added superscript $\delta$ to symbolise its correspondence to a discontinuous solution. An approximation to the continuous flux gradient is then formed via the correction procedure as
    \begin{equation}\label{eq:cont_flux}
        \px{\hat{f}}{\xi} = \px{\hat{f}^\delta}{\xi} + (\hat{f}^\mathrm{num}_L - \hat{f}^\delta_L)\dx{h_L}{\xi} + (\hat{f}^\mathrm{num}_R - \hat{f}^\delta_R)\dx{h_R}{\xi}.
    \end{equation}
    The last and penultimate terms on the right-hand side are used to correct the discontinuous flux to continuous; required for the method to be conservative. Subscripts $L$ and $R$ are used to denote a quantity at the left or right interface and in the case of $\hat{f}^\delta$ these values are interpolated. For $\hat{f}^\mathrm{num}$ these are interface values that are common to the all the elements that share that interface point. Later for linear equations we will define how this is set, but for alternative equation sets approximate Riemann solver offer a suitable means of setting $\hat{f}^\mathrm{num}$. 
    
    In \cref{eq:cont_flux} we introduced the functions $h_L$ and $h_R$, these are the correction functions with the boundary conditions $h_L(-1)=h_R(1)=1$ and $h_R(-1)=h_L(1)=0$. Due to $u_i\in\mathbb{P}_k$, we have that $h_L,h_R\in\mathbb{P}_{k+1}$. Although it is not necessary, it is typical to set $h_L$ and $h_R$ as degree $k+1$ polynomials. The primary aim of this paper is then how to set these correction functions such that methods are linearly stable.
    
    Lastly, the corrected flux gradient can be coupled to an explicit time integration method, such as SSP-RK3, or coupled to some more complex implicit time-integration system, such as diagonally implicit Runge--Kutta~\citep{Wang2020}.

\subsection{Correction Functions}
    In the earliest paper on the subject of FR, \citet{Huynh2007} noted the different numerical properties realised by changing the correction function. In the later works of Vincent \etal~\citep{Vincent2010,Vincent2011}, these correction functions were extended to form continuous classes of functions in one-dimension. With all but one of the correction functions put forward by \citet{Huynh2007} being found to be in that class.
    In the more recent work by \citet{Trojak2021}, a weighted norm was used in the continuous analysis framework of \citet{Vincent2010} to produce yet another one parameter family of correction functions with Jacobi orthogonal polynomials as the basis. This idea can be taken to the natural conclusion for any weight function that is positive almost every by using the three-term recurrence relation to generate sets of monic orthogonal polynomials~\citep{Gautschi1968}.
    
    An alternative approach that was taken by \citet{Trojak2019} was to extend the norm used to define stability. Previously, a limited Sobolev norm was used that was sufficient to define the topology of the approximation space, but does not fully capture it. The high order terms at the interfaces that occur in the stability analysis can not be reconciled with the analysis of \citet{Vincent2011}, it nonetheless showed that it was possible to construct vast sets of correction functions.
    
    In the work of \citet{Castonguay2011}, an analogous family of correction functions to \citet{Vincent2010} was defined on equilateral triangles. To define this, first consider the $(m,k)$-differentiation operator in two dimensions as
    \begin{equation}
        D^{m,k} = \frac{\partial^k}{\partial x^{p-m+1}\partial y^{m-1}}.
    \end{equation}
    The Dubiner basis~\citep{Dubiner1991} can then be defined as
    \begin{definition}[Dubiner Basis]\label{def:dubiner}
        The set
        \begin{equation}
            \phi_i(x,y) = \frac{2}{3^{1/4}}(1-b)^v\psi^{(0,0)}_v(a)\psi_w^{(2v+1,0)}(b)
        \end{equation}
        is orthogonal over a reference equilateral triangle, where 
        \begin{equation}
            a = \frac{3x}{2-\sqrt{3}y}, \quad \mathrm{and}\quad b = \frac{2\sqrt{3}y-1}{3}.
        \end{equation}
        The orders $w$ and $v$, are then integer solutions to:
        \begin{equation}
            0= v^2 - (2k+3)v - 2(w-i+1)\quad\mathrm{for}\quad 0\leq v,w \quad \mathrm{and} \quad w+v\leq k.
        \end{equation}
        where $\psi^{(\alpha,\beta)}_i$ is a normalised Jacobi polynomial and $\psi_i^{(0,0)}=\psi_i$ is a normalised Legendre polynomial. With this definition we can define the set of basis  polynomial as
        \begin{equation}
            \mathbb{Q}_k = \{\phi_i\}_{i\leq(k+1)(k+2)/2}
        \end{equation}
    \end{definition}
    
    We can now define the correction function family of \citet{Castonguay2011} in the following lemma
    \begin{lemma}[Castonguay \etal correction functions]
        For a flux point j at the $\ob{x}_j$ and with surface quadrature wight $w_j$, then defining the reconstructed divergence of the correction corresponding to point $j$ as
        \begin{equation}
            \nabla\cdot\mathbf{h}_j = \sum^N_{i=1}\sigma_i\phi_i,
        \end{equation}
        then if the modal coefficients are found from
        \begin{equation}
            c\sum^{N}_{l=1}\sigma_l\sum^{k+1}_{m=1}\binom{k}{m-1}D^{m,k}\phi_iD^{m,k}\phi_l=-\sigma_i+w_j\delta_{ij}\phi_i(\ob{x}_j),
        \end{equation}
        a sufficient condition for linear stability of the resulting FR method is that $c>0$.
    \end{lemma}
    \begin{proof}
        See \citet{Castonguay2011}.
    \end{proof}
    \begin{remark}
        Both the Dubiner basis and the Castonguay correction functions are defined on a triangle with a total order basis~\citep{Trefethen2017}. This is the most commonly used basis for simplex elements due to the trace space being a polynomial space, and so we will restrict our investigation to these elements.
    \end{remark}
    
    The work discussed until now used a continuous approach to find correction functions. A different insight may be gained if a discrete approach is used. This was the approach used by \citet{Vincent2015} to produce an extended range of stable corrections, which was later encompassed in the work of \citet{Trojak2019}. This discrete approach was further formalised within the summation-by-parts framework in the work of \citet{Ranocha2016}. In that work it was also shown how a skew-symmetric split form with a lumped mass matrix could be used to prove stability for Burgers' equation, but not in the general case.
    
\subsection{Summation-by-parts}
    Many advances have been made in the theory of DG and FR methods by considering the discrete problem. A foundation of these analytic techniques is the definition of summation-by-parts (SBP) operators.
    
    First let us define the basic discrete operators that will be used throughout this work. If $u_i$ is an approximation in element $K_i$ to the function $u\in C^1(K)$, where for our domain we have $K\subset \mathbb{R}^d$. Then for some set of $N_s$ solution points $\mathbf{x}_i=\{\mathbf{x}_{i,j}\}_{j\leq N_s}$ we have the vector $\mathbf{u}_i=u_i(\mathbf{x}_i)$, which is the evaluation of $u_i$ at the solution points. If $l_j(\mathbf{x})$ is the Lagrange polynomial corresponding to $l_j(\mathbf{x}_{i,k})=\delta_{jk}$, then we can define a mass matrix as
    \begin{equation}\label{eq:mass_matrix}
        \mathbf{u}_i^T\mathbf{M}\mathbf{u}_i = \int_{K_i} \left(\sum^{N_s}_{j=1}u_i(\mathbf{x}_{i,j})l_j(\mathbf{x})\right)^2\mathrm{d}\mathbf{x}.
    \end{equation}
    If we call the cardinal axes $x_1,\dots,x_d$, then we can define the differentiation operators
    \begin{equation}
        \mathbf{D}_{x_1}\mathbf{u}_i = \sum^{N_s}_{j=1}u_i(\mathbf{x}_{i,j})\dx{l_j(\mathbf{x})}{x_1}, \quad \mathbf{D}_{x_2}\mathbf{u}_i = \sum^{N_s}_{j=1}u_i(\mathbf{x}_{i,j})\dx{l_j(\mathbf{x})}{x_2}, \quad \dots
    \end{equation}
    In the following we will drop the subscript $i$ for clarity except where it is explicitly needed.
    
    Moving on to define SBP in higher-dimensions, we start with the analogy of integration-by-parts in higher dimensions
    \begin{definition}[Divergence integration-by-parts]
        For a scalar field, $v\in C^1(K)$ and a vector field, $\mathbf{W}\in (C^1(K))^d$, in the closed domain $K\subset\mathbb{R}^d$ with boundary $\partial K$, then
        \begin{equation}
            \int_K v\nabla\cdot\mathbf{W}\mathrm{d}K + \int_K\nabla v\cdot\mathbf{W}\mathrm{d}K = \int_{\partial K} v\mathbf{W}\cdot \mathbf{n} \mathrm{d}S.
        \end{equation}
    \end{definition}
    With this we may then define the generalised SBP relation as
    \begin{definition}[Generalised Summation-by-parts]\label{def:sbp}
        For solutions $u\in C^1(K)$ and $U\in (C^1(K))^d$, let  $u_i$ and $U_i$ be approximations in element $K_i$ such that for some nodal point set $\mathbf{x}_i\{\mathbf{x}_{i,j}\}_{j\leq N_s}$ we have $\mathbf{u}_i = u_i(\mathbf{x}_i)$ and $\mathbf{U}_i = U_i(\mathbf{x}_i)$, then a set of operators is said to satisfy the generalised SBP property if:
        \begin{equation}\label{eq:sbp}
            \mathbf{MD} + \mathbf{G}^T\hb{M} = \mathbf{L}^T_{\partial K}\mathbf{W}_{\partial K}\mathbf{N}\hb{L}_{\partial K},
        \end{equation}
        where the divergence and gradient operators are defined as
        \begin{equation}
            \mathbf{DU}_i = [\mathbf{D}_{x_1}, \mathbf{D}_{x_2}, \dots]\mathbf{U}_i = \nabla\cdot U_i \quad \mathrm{and} \quad \mathbf{Gu}_i = 
            \begin{bmatrix}
                \mathbf{D}_{x_1}\\ \mathbf{D}_{x_2}\\\vdots
            \end{bmatrix}\mathbf{u} = \nabla u_i.
        \end{equation}
        Then defining the interpolation operator $\mathbf{L}_{\partial K}:K\mapsto\partial K$, and surface mass matrix, $\mathbf{W}_{\partial K}$, such that
        \begin{equation}
            \mathbf{u}_i^T\mathbf{L}_{\partial K}\mathbf{W}_{\partial K}\mathbf{N}\hb{L}_{\partial K}\mathbf{U}_i = \int_{\partial K}u_iU_i\cdot \mathbf{n}_i\mathrm{d}s,
        \end{equation}
        where $\mathbf{n}$ is a vector of outwards facing normals at the interface. Finally, the Kronecker product of a matrix with the identity matrix is denoted by
        \begin{equation}
            \hb{A} = \mathbf{I}_d\otimes\mathbf{A}.
        \end{equation}
    \end{definition}
    \begin{remark}
        From this definition of SBP we see that the restriction on the mass matrix is that it should accurately integrate all functions in at least $\mathbb{Q}_{2k-1}$. From the definition of the mass matrix in \cref{eq:mass_ma} this is true, however, in many applications this may not be true if using a quadrature. This case is explicitly handled in \cref{app:quad}.
    \end{remark}
    
    SBP is simply a discrete restatement of integration-by-parts. The advantage is it permits the development of discrete analogues of continuous properties of the physical system. The earliest discussion of SBP in the context of finite element methods, to the authors' knowledge, is that of \citet{Fisher2013}. This was an adaptation of ideas previously used throughout the finite difference community. There are many works studying SBP in a finite difference context, with some of the earliest works including \citet{Carpenter1995} and \citet{Olsson1995}. An important work in concreting the utility of these approaches is that of \citet{Carpenter1994}. There it was shown that for finite differences applied to hyperbolic systems, a scheme with energy bounded via SBP leads to the solution being bounded in the continuous problem. This is important as it shows consistency of the discrete stability analysis and the continuous problem.

    In general the exploration of SBP operators has largely focused on one-dimension, but some recent works have move beyond this. For example, on tensor-product elements~\citep{DelReyFernandez2019}. In the work of \citet{Hicken2016} they were able to extend the theory to simplex elements using the generalisation of SBP by \citet{DelReyFernandez2014}. The form given in \citep{DelReyFernandez2014} is analogous to that shown in \cref{def:sbp}.
    
    The operators set out in \cref{def:sbp} can then be used to construct the FR scheme. First consider a linear advection equation such as
    \begin{equation}\label{eq:cons_eq}
        \px{u}{t} + \nabla\cdot\mathbf{f} = 0,
    \end{equation}
    where $\mathbf{f}=\mathbf{a}u$. The FR discretisation of the spatial derivatives can be written within the SBP framework as:
    \begin{equation}
        \px{\mathbf{u}}{t} = -\mathbf{D}\mathbf{F} - \mathbf{C}\left[ (\mathbf{n}\cdot\mathbf{F}^\mathrm{num} - \mathbf{N}\hb{L}_\partial\mathbf{F})\right]
    \end{equation}
    where $\mathbf{C}$ is the correction function matrix.
    
    In the previous work of \citet{Castonguay2011} and \citet{Vincent2015}, the modal form was used in the presentation of the stable correction functions as the forms are far sparser. Transformation from a nodal to modal representation is defined via
    \begin{equation}
        \mathbf{u} = \mathbf{V}\tb{u},
    \end{equation}
    where $\tb{u}$ is a vector of modal coefficients, and $\mathbf{V}$ is the Vandermonde matrix. Throughout the rest of this work we will use a tilde to denote a matrix of vector in the \emph{modal} representation.
\section{Linear stability analysis}\label{sec:linear}
    To study linear stability, we prescribe the system being solved as a generalised linear advection equation with the form:
    \begin{equation}\label{eq:lin_advec}
        \px{u}{t} + \nabla\cdot(u\otimes \mathbf{a}) = 0, \quad \mathrm{for} \quad  \mathbf{a} = [a_{x_1},a_{x_2},\dots]^T.
    \end{equation}
    
    First we consider the known correction function where FR corresponds to DG, here stability is found in the $L_2$ norm induced by $\mathbf{M}$. We do this to demonstrate the use of the SBP framework in higher dimensions and to more clearly set out the interface treatment.
    \begin{lemma}[Linear Stability]\label{lem:lin_stab}
        Setting the nodal solution as $\mathbf{u}_i=u(\mathbf{x}_i)$ and linear nodal flux as $\mathbf{F} = [a_{x_1},a_{x_2},\dots]^T\otimes\mathbf{u}$, then for the FR scheme applied to \cref{eq:lin_advec}, the following conditions
        \begin{equation}\label{eq:la_cond1}
            \mathbf{C} = \mathbf{M}^{-1}\mathbf{L}_\partial^T\mathbf{W}_\partial,
        \end{equation}
        \begin{equation}\label{eq:la_cond2}
            \mathbf{W}_\partial = \mathrm{diag}(\mathbf{w}_\partial), \quad \text{and} \quad \mathbf{w}_i\geq 0,
        \end{equation}
        and
        \begin{subequations}\label{eq:la_cond3}
            \begin{align}
                (\mathbf{n}\cdot F)^{\mathrm{num} +}_j &= \half(\mathbf{n}_j^+\cdot\mathbf{a})(u_j^+ + u_j^-) - \half\kappa|\mathbf{n}_j^+\cdot\mathbf{a}|(u_j^- - u_j^+), \quad \text{and}\\
                (\mathbf{n}\cdot F)^{\mathrm{num} -}_j &= \half(\mathbf{n}_j^-\cdot\mathbf{a})(u_j^- + u_j^+) - \half\kappa|\mathbf{n}_j^-\cdot\mathbf{a}|(u_j^+ - u_j^-), \quad \text{for} \quad \kappa \in [0,1],
            \end{align}
        \end{subequations}
        are sufficient for energy stability in the norm induced by the mass matrix, $\mathbf{M}$, \ie
        \begin{equation}
            \dx{}{t}\|\mathbf{u}\|_M^2 \leq 0.
        \end{equation}
    \end{lemma}
    \begin{proof}
        The FR method applied to \cref{eq:lin_advec} can be written as
        \begin{equation*}
            \dx{\mathbf{u}}{t} = -\mathbf{DF} - \mathbf{C}\left[(\mathbf{n}\cdot \mathbf{F})^\mathrm{num} - \mathbf{N}\hb{L}_\partial\mathbf{F}\right].
        \end{equation*}
        Then multiplying this by $\mathbf{u}^T\mathbf{M}$ we get
        \begin{equation}\label{eq:la_eq1}
            \mathbf{u}^T\dx{}{t}\mathbf{Mu} = -\mathbf{u}^T\mathbf{MDF} - \mathbf{u}^T\mathbf{MC}\left[(\mathbf{n}\cdot \mathbf{F})^\mathrm{num} - \mathbf{N}\hb{L}_\partial\mathbf{F}\right],
        \end{equation}
        we can then apply \cref{eq:sbp} to obtain a second equation
        \begin{equation}\label{eq:la_eq2}
            \mathbf{u}^T\dx{}{t}\mathbf{Mu} = \mathbf{u}^T\mathbf{G}^T\hb{M}\mathbf{F} - \mathbf{u}^T\mathbf{L}^T_\partial\mathbf{W}_\partial\mathbf{N}\hb{L}_\partial\mathbf{F} - \mathbf{u}^T\mathbf{MC}\left[(\mathbf{n}\cdot \mathbf{F})^\mathrm{num} - \mathbf{N}\hb{L}_\partial\mathbf{F}\right].
        \end{equation}
        \cref{eq:la_eq1,eq:la_eq2} can then be combined to give
        \begin{equation}
            2\mathbf{u}^T\dx{}{t}\mathbf{Mu} = \dx{}{t}\|\mathbf{u}\|_M^2 =  - \mathbf{u}^T\mathbf{L}^T_\partial\mathbf{W}_\partial\mathbf{N}\hb{L}_\partial\mathbf{F} - 2\mathbf{u}^T\mathbf{MC}\left[(\mathbf{n}\cdot \mathbf{F})^\mathrm{num} - \mathbf{N}\hb{L}_\partial\mathbf{F}\right]
        \end{equation}
        here we have used the symmetry of $\mathbf{M}$ which leads to $\mathbf{u}^T\mathbf{G}^T\hb{M}\mathbf{F} - \mathbf{u}^T\mathbf{MDF}=0$. Then applying \cref{eq:la_cond1} we recover
        \begin{equation}\label{eq:la_eq3}
            \dx{}{t}\|\mathbf{u}\|_M^2 = \mathbf{u}^T\mathbf{L}^T_\partial\mathbf{W}_\partial\mathbf{N}\hb{L}_\partial\mathbf{F} - 2\mathbf{u}^T\mathbf{L}^T_\partial\mathbf{W}_\partial(\mathbf{n}\cdot\mathbf{F})^\mathrm{num} = \mathbf{u}^T\mathbf{L}^T_\partial\mathbf{W}_\partial(\mathbf{N}\hb{L}_\partial\mathbf{F} - 2(\mathbf{n}\cdot\mathbf{F})^\mathrm{num}).
        \end{equation}
        Now considering a mesh of multiple elements and focusing on a single point on the boundary of an element, say point $j$. Using the condition that $\mathbf{W}_\partial$ is diagonal, then the global contribution to the right-hand-side of  \cref{eq:la_eq3} from point $j$ is:
        \begin{equation}
            u_j^-w_j^-\left[(\mathbf{n}_j^-\cdot F_j^-) - 2(\mathbf{n}\cdot F)^{\mathrm{num}-}_j\right] +  u_j^+w_j^+\left[(\mathbf{n}^+_j\cdot F_j^+) - 2(\mathbf{n}\cdot F)^{\mathrm{num}+}_j\right].
        \end{equation}
        Here $-$ and $+$ are the contributions from either side of the interface, and from \cref{eq:la_cond2} $w_j$ is the positive surface quadrature weight at $j$. Then setting the numerical flux from \cref{eq:la_cond3} we obtain
        \begin{multline}
            u_j^-w_j^-\left[(\mathbf{n}_j^-\cdot F_j^-) - 2(\mathbf{n}\cdot F)^{\mathrm{num}-}_j\right] +  u_j^+w_j^+\left[(\mathbf{n}^+_j\cdot F_j^+) - 2(\mathbf{n}\cdot F)^{\mathrm{num}+}_j\right] \\
            = u_j^+u_j^-(\mathbf{n}_j\cdot\mathbf{a})(w_j^--w_j^+) - \kappa|\mathbf{n}_j\cdot\mathbf{a}|(u_j^--u_j^+)(w_j^-u_j^--w_j^+u_j^+),
        \end{multline}
        where we have used $\mathbf{n}_j^+=-\mathbf{n}_j^-$ by definition for a conformal mesh. Applying \cref{eq:la_cond4} we recover
        \begin{equation}
            u_j^+u_j^-(\mathbf{n}_j\cdot\mathbf{a})(w_j^--w_j^+) - \kappa|\mathbf{n}_j\cdot\mathbf{a}|(u_j^--u_j^+)(w_j^-u_j^--w_j^+u_j^+) = -\kappa|\mathbf{n}_j\cdot\mathbf{a}|w^+_j(u_j^--u_j^+)^2\leq0,
        \end{equation}
        where we further assume that $w_j^-=w_j^+$ \ie the flux points used here have some degree of rotational symmetry.
        And hence summing over the domain we recover the required result of
        \begin{equation}
            \dx{}{t}\|\mathbf{u}\|^2_M \leq 0.
        \end{equation}
    \end{proof}
    
    As described in \cref{sec:prelim}, \citet{Vincent2015} and later  \citet{Ranocha2016} were able to derive a multi-parameter extended-range set of FR methods in one dimension that were linearly stable. These methods were found to be stable in a modified norm such that:
    \begin{equation}
        \dx{}{t}\mathbf{u}^T(\mathbf{M}+\mathbf{Q})\mathbf{u} = \dx{}{t}\|\mathbf{u}\|_{M+Q}^2 \leq 0
    \end{equation}
    We now generalise this set of methods to higher dimensions in the following lemma
    \begin{lemma}[Extended-range linear stability]\label{lem:lin_stab_q}
        For the conditions set out in \cref{lem:lin_stab}, with the additional constraints that 
        \begin{subequations}\label{eq:la_cond_q}
            \begin{align}
                \mathbf{Q} - \mathbf{Q}^T & = 0\label{eq:la_cond4}\\
                \mathbf{QD}_{x_i} + \mathbf{D}^T_{x_i}\mathbf{Q}^T &= 0 \quad \forall \; i\in\{1,\cdots,d\}\label{eq:la_cond5}\\
                \mathbf{v}^T(\mathbf{M} + \mathbf{Q})\mathbf{v} &> 0, \quad \forall \; \mathbf{v}\in\mathbb{R}^{N_k}\setminus 0\label{eq:la_cond6}
            \end{align}
        \end{subequations}
        and the modified condition that
        \begin{equation}\label{eq:la_cond7}
            \mathbf{C} = (\mathbf{M} + \mathbf{Q})^{-1}\mathbf{L}_\partial^T\mathbf{W}_\partial,
        \end{equation}
        then FR applied to \cref{eq:lin_advec} is stable in the norm induced by $(\mathbf{M}+\mathbf{Q})$, \ie
        \begin{equation}
            \dx{}{t}\|\mathbf{u}\|^2_{M+Q} \leq 0.
        \end{equation}
    \end{lemma}
    \begin{proof}
        Following the same steps as in the proof of \cref{lem:lin_stab} and using the modified condition in \cref{eq:la_cond7} we obtain
        \begin{equation}
            \dx{}{t}\|\mathbf{u}\|_{M+Q}^2 = -2\mathbf{u}^T\mathbf{QDF} - \mathbf{u}^T\mathbf{L}^T_\partial\mathbf{W}_\partial\mathbf{N}\hb{L}_\partial\mathbf{F} - 2\mathbf{u}^T(\mathbf{M}+\mathbf{Q})\mathbf{C}\left[(\mathbf{n}\cdot \mathbf{F})^\mathrm{num} - \mathbf{N}\hb{L}_\partial\mathbf{F}\right]
        \end{equation}
        Then applying \cref{eq:la_cond4,eq:la_cond5}, the first term of the right-hand-side is found to be zero, and so proceeding with the proof of \cref{lem:lin_stab}, we achieve the required results of 
        \begin{equation}
            \dx{}{t}\|\mathbf{u}\|^2_{M+Q} \leq 0
        \end{equation}
        The final condition \cref{eq:la_cond6} is used to ensure that the norm induced by $\mathbf{M}+\mathbf{Q}$ is a valid norm.
    \end{proof}
    \begin{remark}
        A stricter condition on $\mathbf{Q}$ is $\mathbf{QD} =- \mathbf{G}^T\hb{Q}$; however, when looking for a $\mathbf{Q}$ that satisfies this the solution $\mathbf{Q}=\mathbf{0}$ is typically recovered. Alternatively, if the less strict condition \cref{eq:la_cond5} is used, a wider range of valid $\mathbf{Q}$ matrices can be found that still guarantee linear stability.
    \end{remark}
    \begin{remark}
        By finding a norm $\|\mathbf{u}\|_{M+Q}$ where the solution norm monotonically decays in time, we can use the equivalence of norms to infer stability. Therefore, as $c\|\mathbf{u}\|_M\leq\|\mathbf{u}|_{M+Q}\leq C\|\mathbf{u}\|_M$, the norm $\|\mathbf{u}\|_M$ may not decay monotonically in time, but its rate of decay must remain bounded.
    \end{remark}

    It is often convenient to consider methods in the modal form rather than the nodal form, but to be confident that a scheme found to be stable in modal form is stable in nodal form consider the following:
    \begin{corollary}[Nodal-modal equivalence]\label{col:nodal_modal}
        The stability of a scheme that satisfies conditions \cref{eq:la_cond2,eq:la_cond3,eq:la_cond_q,eq:la_cond7} is independent of modal or nodal representation, provided $\mathbf{V}$ is invertible.
    \end{corollary}
    \begin{proof}
        To prove this it is sufficient to show that the conditions \cref{eq:la_cond_q}, if satisfied in one frame, are satisfied in the other. First consider the transform of $\mathbf{Q}$ as
        \begin{equation}
            \tb{Q} = \mathbf{V}^T\mathbf{Q}\mathbf{V},
        \end{equation}
        clearly if $\mathbf{Q}=\mathbf{Q}^T$ then $\tb{Q}=\tb{Q}^T$. Next considering the skew symmetry property we have
        \begin{equation}
            \mathbf{QD} = \mathbf{V}^{-T}\tb{Q}\mathbf{V}^{-1}\mathbf{V}\tb{D}\mathbf{V}^{-1} = \mathbf{V}^{-T}\tb{Q}\tb{D}\mathbf{V}^{-1}
        \end{equation}
        and 
        \begin{equation}
            -\mathbf{D}^T\mathbf{Q} = -(\mathbf{V}\tb{D}\mathbf{V}^{-1})^T\mathbf{V}^{-T}\tb{Q}\mathbf{V}^{-1} = -\mathbf{V}^{-T}\tb{D}^T\tb{Q}\mathbf{V}^{-1}.
        \end{equation}
        Therefore, if $\mathbf{QD} = -\mathbf{D}^T\mathbf{Q}$, then $\tb{Q}\tb{D}=-\tb{D}^T\tb{Q}$. Lastly for \cref{eq:la_cond6} we have to show that if $\mathbf{M}+\mathbf{Q}$ is positive definite, then so is $\tb{M}+\tb{Q}$. Considering this property we have that:
        \begin{equation}
            \mathbf{w}^T(\mathbf{M}+\mathbf{Q})\mathbf{V}^{-1}\mathbf{w} = \mathbf{w}^T\mathbf{V}^{-T}(\tb{M}+\tb{Q})\mathbf{V}^{-1}\mathbf{w} = \mathbf{v}^T(\tb{M}+\tb{Q})\mathbf{v} > 0,
        \end{equation}
        which holds as $\mathbf{V}$ is full rank. This completes the proof.
    \end{proof}

\subsection{Conservation}
    \begin{lemma}[Linear conservation]\label{lem:lin_cons}
        Consider the $d$-dimensional FR method with linear flux function such that, for Banach space $V$, $\mathbf{u}\in V$ and $\mathbf{F} = F(\mathbf{u})\in (V^\prime,V)^d$, then sufficient conditions for conservation are that the gradient operator is are least first order accurate, \ie
        \begin{equation}\label{eq:sc_cond1}
            \mathbf{G1} = \mathbf{0},
        \end{equation}
        and that the lifting operator is such that
        \begin{equation}\label{eq:sc_cond2}
            \mathbf{1}^T\mathbf{MC} = \mathbf{1}^T\mathbf{L}^T_\partial\mathbf{W}_\partial.
        \end{equation}
    \end{lemma}
    \begin{proof}
        Let $\mathbf{u}_i=u(\mathbf{x}_i)$ and $\mathbf{F} = F(\mathbf{u})$, then the FR method can be written as
        \begin{equation}
            \dx{\mathbf{u}}{t} = -\mathbf{DF} - \mathbf{C}\left[(\mathbf{n}\cdot\mathbf{F})^\mathrm{num} - \mathbf{N}\hb{L}_{\partial}\mathbf{F}\right]
        \end{equation}
        Then multiplying by $\mathbf{1}^T\mathbf{M}$ we obtain
        \begin{equation*}
            \mathbf{1}^T\dx{}{t}\mathbf{Mu}= \dx{}{t}\mathbf{1}^T\mathbf{Mu} = -\mathbf{1}^T\mathbf{MDF} - \mathbf{1}^T\mathbf{MC}\left[(\mathbf{n}\cdot\mathbf{F})^\mathrm{num} - \mathbf{N}\hb{L}_{\partial}\mathbf{F}\right].
        \end{equation*}
        Then applying \cref{eq:sbp} we obtain
        \begin{equation*}
            \dx{}{t}\mathbf{1}^T\mathbf{Mu} = \mathbf{1}^T\mathbf{G}^T\hb{M}\mathbf{F} - \mathbf{1}^T\mathbf{L}_\partial^T\mathbf{W}_\partial\mathbf{N}\hb{L}_\partial\mathbf{F} - \mathbf{1}^T\mathbf{MC}\left[(\mathbf{n}\cdot\mathbf{F})^\mathrm{num} - \mathbf{N}\hb{L}_{\partial}\mathbf{F}\right].
        \end{equation*}
        If \cref{eq:sc_cond1} holds, then we obtain
        \begin{equation*}
            \dx{}{t}\mathbf{1}^T\mathbf{Mu} = - \mathbf{1}^T\mathbf{L}_\partial^T\mathbf{W}_\partial\mathbf{N}\hb{L}_\partial\mathbf{F} - \mathbf{1}^T\mathbf{MC}\left[(\mathbf{n}\cdot\mathbf{F})^\mathrm{num} - \mathbf{N}\hb{L}_{\partial}\mathbf{F}\right],
        \end{equation*}
        and proceeding to apply \cref{eq:sc_cond2} we get
        \begin{equation*}
            \dx{}{t}\mathbf{1}^T\mathbf{Mu} = - \mathbf{1}^T\mathbf{L}^T_\partial\mathbf{W}_\partial(\mathbf{n}\cdot\mathbf{F})^\mathrm{num}.
        \end{equation*}
        The term on the right-hand side is discrete statement of divergence theorem, which in 1D would give $f_R - f_L$. Therefore, the scheme is conservative. 
    \end{proof}
    \begin{remark}
        A similar lemma can be defined for non-linear flux functions, if an intermediate set of quadrature points is used, see \citet{Chan2018}.
    \end{remark}
    
    Here we set $\mathbf{Q}=0$, but \citet{Vincent2015} showed how changing $\mathbf{Q}$ could lead some methods to be non-conservative in an integral with unit measure. Conservation of the extended range of stable schemes is then considered in the following lemma
    \begin{lemma}[Conservation of extended schemes]
        For an FR scheme that satisfies \cref{eq:la_cond2,eq:la_cond3,eq:la_cond_q,eq:la_cond7}, with a linear flux function, then if the following condition is also satisfied
        \begin{equation}
            \mathbf{1}^T\mathbf{M}(\mathbf{M}+\mathbf{Q})^{-1}\mathbf{L}^T_\partial\mathbf{W}_\partial = \mathbf{1}^T\mathbf{L}^T_\partial\mathbf{W}_\partial
        \end{equation}
        then the scheme is conservative in that
        \begin{equation}
            \dx{}{t}\mathbf{1}^T\mathbf{Mu} = \mathbf{f}_l - \mathbf{f}_r
        \end{equation}
    \end{lemma}
    This can be straightforwardly seen from \cref{lem:lin_cons} and \cref{eq:la_cond7}. 
    \begin{remark}
        What can be seen from \cref{lem:lin_cons} and \cref{lem:lin_stab_q} is that any FR method satisfying \cref{eq:la_cond2,eq:la_cond3,eq:la_cond_q,eq:la_cond7} for a linear flux, is automatically conservative in terms of the norm induced by $\mathbf{M}+\mathbf{Q}$. However, this is not physical and could lead to schemes that are not consistent.
    \end{remark}
    
\subsection{Symmetry Conditions}
    It is taken as axiomatic that the numerical method should be independent of node ordering, or element orientation. As the  correction function can change the numerical characteristics of the FR method; therefore, the correction function is required to have some degree of symmetry.
    
\begin{figure}[tbhp]
        \centering
        \subfloat[Axis rotation.]{\label{fig:syms}\adjustbox{width=0.25\linewidth,valign=b}{\begin{tikzpicture}[scale=2]
    
    \begin{scope}[on behind layer]
        \draw[black, -, thick] (-1,{-1/sqrt(3)}) --(1,{-1/sqrt(3)}) -- (0,{2/sqrt(3)}) -- cycle;
    \end{scope}
    
    \begin{scope}
        \draw[fill={Pastel1-A}, color={Pastel1-A}] (-1,{-1/sqrt(3)}) circle (1.5pt);
        \draw[fill={Pastel1-A}, color={Pastel1-A}] (1,{-1/sqrt(3)}) circle (1.5pt);
        \draw[fill={Pastel1-A}, color={Pastel1-A}] (0,{2/sqrt(3)}) circle (1.5pt);
    \end{scope}
    
    \draw[black, -latex] (0,-0.52735) -- (0,-0.27735) node[left] {\footnotesize{$y$}};
    \draw[black, -latex] (0,-0.52735) -- (0.25,-0.52735) node[right] {\footnotesize{$x$}};
    
    \draw[black, -latex] (0.4567,0.26368) -- (0.24019,0.13868) node[below,shift={(0.15,0.05)}] {\footnotesize{$y^{\prime}$}};
    \draw[black, -latex] (0.4567,0.26368) -- (0.3317,0.48018) node[above,shift={(-0.17,-0.06)}] {\footnotesize{$x^{\prime}$}};
    
    \draw[black,-] (0,-0.2) -- (0,0) -- ({0.2*0.5*sqrt(3)},{0.5*0.2});
    \draw[black, thick, densely dotted] (0,-0.17) arc (-90:30:0.17) node[midway, shift={(0.1,-0.1)}] {\footnotesize{$\theta$}};
    
     \begin{scope}[on behind layer]
        \draw[black!00] (-1,-0.775) rectangle (1,1.1547);
    \end{scope}
\end{tikzpicture}}}
        \subfloat[Node locations and rotational symmetries.]{\label{fig:ref_triangle}\adjustbox{width=0.45\linewidth,valign=b}{\begin{tikzpicture}[scale=2]

    \tikzstyle{xnode}=[draw, circle, thick, color={Pastel1-A}, fill={Pastel1-A}, text=black, minimum width=1pt]
                      
    \tikzstyle{edge}=[-,black,thick]
    
    \begin{scope}
        \draw[fill={Pastel1-A}, color={Pastel1-A}] (-1,-0.57735026919) circle (1.5pt) node [left, text=black, shift={(-0.1,0.0)}] {\footnotesize{$\left(-1,-\frac{1}{\sqrt{3}}\right)$}};
        
        \draw[fill={Pastel1-A}, color={Pastel1-A}] ( 1,-0.57735026919) circle (1.5pt) node [right, text=black, shift={( 0.1,0.0)}] {\footnotesize{$\left(1,-\frac{1}{\sqrt{3}}\right)$}};
        
        \draw[fill={Pastel1-A}, color={Pastel1-A}] ( 0, 1.15470053838) circle (1.5pt) node [above, text=black, shift={( 0.0,0.1)}] {\footnotesize{$\left(0,\frac{2}{\sqrt{3}}\right)$}};
    \end{scope}
    
    \begin{scope}[on behind layer]
        \draw[thick, black] (-1,-0.57735026919) -- ( 1,-0.57735026919) -- ( 0, 1.15470053838) -- cycle;
    \end{scope}
    
    
    \draw[black, -latex] (0,-0.5) -- (0,-0.2) node[left] {\footnotesize{$y$}};
    \draw[black, -latex] (0,-0.5) -- (0.3,-0.5) node[right] {\footnotesize{$x$}};
    
    \draw[black, -latex] (-0.45669872981,0.26367513459) -- (-0.19689110867,0.11367513459) node[above, shift={( 0.0,0.1)}] {\footnotesize{$y^\prime$}};
    \draw[black, -latex] (-0.45669872981,0.26367513459) -- (-0.60669872981,0.00386751345) node[below,, shift={( 0.0,0.1)}] {\footnotesize{$x^\prime$}};
    
    \draw[black, -latex] ( 0.45669872981,0.26367513459) -- ( 0.19689110867,0.11367513459) node[below, shift={( 0.2,0.1)}] {\footnotesize{$y^{\prime\prime}$}};
    \draw[black, -latex] ( 0.45669872981,0.26367513459) -- ( 0.30669872981,0.52348275573) node[left, shift={( 0.1,0.05)}] {\footnotesize{$x^{\prime\prime}$}};
    
\end{tikzpicture}}}
        ~
        \subfloat[Symmetry axes.]{\label{fig:triangle_axes}\adjustbox{width=0.25\linewidth,valign=b}{\begin{tikzpicture}[scale=2]

    \tikzstyle{xnode}=[draw, circle, thick, color={Pastel1-A}, fill={Pastel1-A}, text=black, minimum width=1pt]
                      
    \tikzstyle{edge}=[-,black,thick]
    
    \begin{scope}
        \draw[fill={Pastel1-A}, color={Pastel1-A}] (-1,-0.57735026919) circle (1.5pt);
        
        \draw[fill={Pastel1-A}, color={Pastel1-A}] ( 1,-0.57735026919) circle (1.5pt);
        
        \draw[fill={Pastel1-A}, color={Pastel1-A}] ( 0, 1.15470053838) circle (1.5pt);
    \end{scope}
    
    \begin{scope}[on behind layer]
        \draw[thick, black] (-1,-0.57735026919) -- ( 1,-0.57735026919) -- ( 0, 1.15470053838) -- cycle;
        \draw[thick, black, densely dotted] (0,-0.75) -- (0,1.3547);
    \end{scope}
    
    \draw[black, -latex] (0.1,-0.5) -- (0.1,-0.2) node[above] {\footnotesize{$y$}};
    \draw[black, -latex] (0.1,-0.5) -- (0.4,-0.5) node[right] {\footnotesize{$x$}};
    
    \draw[black, -latex] (-0.1,-0.5) -- (-0.1,-0.2) node[above] {\footnotesize{$y^\prime$}};
    \draw[black, -latex] (-0.1,-0.5) -- (-0.4,-0.5) node[left] {\footnotesize{$x^\prime$}};
    
    \begin{scope}[on behind layer]
        \draw[black!00] (-1,-0.775) rectangle (1,1.1547);
    \end{scope}
\end{tikzpicture}}}
        \caption{Reference triangle and symmetry definitions.}
    \end{figure}
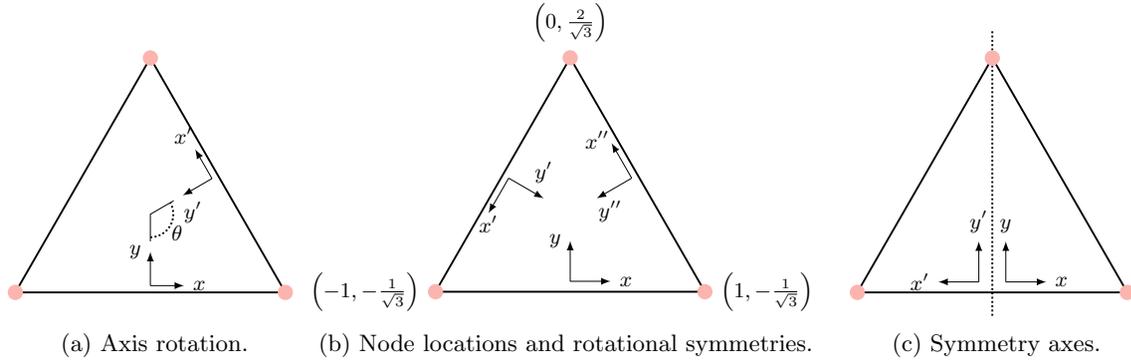

    Defining the cardinal axes for different face reference frames, as in \cref{fig:syms}, the transformation of reference coordinates can be made via:
    \begin{equation}
        x^\prime(\theta) = x\cos\theta + y\sin\theta, \quad \mathrm{and} \quad y^\prime(\theta) = -x\sin\theta + y\cos\theta.
    \end{equation}
    A transformation matrix, $\tb{T}_{mn}$, can then be defined which transforms the basis from the face reference space with $\theta_m$ to $\theta_n$. This allows rotational symmetry conditions to be imposed on $\tb{Q}$ to give
    \begin{equation}
        \tb{T}_{mn}\tb{Q} = \tb{Q}\tb{T}_{mn}.
    \end{equation}
    This condition ensures that a function such as $\phi_i(x,y)$ and the same function rotated to the new reference,  $\phi^\prime_i(x^\prime,y^\prime)$, then have the same value in the norm induced by $\mathbf{M} + \mathbf{Q}$. 
    
    A further symmetry condition is that, given a pair of flux points on a face that are symmetric about some axis, the corresponding correction functions should be symmetric. This gives the condition that 
    \begin{equation}
        \tb{S}_x\tb{Q} = \tb{Q}\tb{S}_x,
    \end{equation}
    where $x$ is the axis of symmetry and $\tb{S}_x$ is a matrix that reflects the modes about the axis $x$. A comparable axial symmetry condition was proposed by \citet{Ranocha2016} for use in one-dimension. Finally, when applying symmetry conditions care must be taken to not over constrain the system. This occurs when one symmetry in the reference frame of a face is a linear combination of other symmetry conditions, and is often indicated by the erroneous recovery of $\tb{Q}:=0$. For the reference triangle this means that applying two rotational symmetries and one axial symmetry is over constrained as one rational symmetry can be expressed using the other rotation and the axial symmetry. 
    
\section{Extended-Range Scheme on Triangle}\label{sec:erfr_tri}
    \citet{Vincent2015} and \citet{Ranocha2016} derived an extended range of energy stable 1D correction functions and the analysis presented in \cref{sec:linear} derived the conditions required to extend this family to triangles. We now set out this generalised extended range of stable correction functions for several orders using the reference triangle shown in \cref{fig:ref_triangle}. Furthermore, schemes will be defined in the modal form due to the sparsity of the matrices, and is supported by \cref{col:nodal_modal}.

    In this section we also look to recover the single parameter schemes of \citet{Castonguay2011}. This set can be cast into the SBP framework with the following definition:
    \begin{definition}[Castonguay \etal simplex method]\label{def:cast_fr}
        Given the reference triangular element of \cref{fig:ref_triangle} and a total order basis, FR is found to be stable in the broken Sobolev norm:
        \begin{equation}
            \half\int_{\hat{K}} u^2 + \frac{c}{A}\sum^{k+1}_{m=1} {{k}\choose{m-1}} \left(D^{m,k}u\right)^2\mathrm{d}\mathbf{x}, 
        \end{equation}
        where $A=|\hat{K}|$. This can then be used to define $\tb{Q}_\mathrm{C}$ required to recover this set of methods in the extended range of stable schemes defined here. Therefore:
        \begin{equation}
            \tb{Q}_\mathrm{C} = \frac{c}{A}\sum^{k+1}_{m=1}{{k}\choose{m-1}}\left(\tb{K}^{m,k}\right)^T\tb{M}\tb{K}^{m,k}, \quad \mathrm{for}\quad \tb{K}^{k,m}=\tb{D}_r^{k-m+1}\tb{D}_s^{m-1}.
        \end{equation}
    \end{definition}
    In what follows we will then look if and how this $\tb{Q}_C$ matrix be recovered in the new set of schemes defined, and what the constraints on $c$ are for stability.
    
\subsection{$k=2$}\label{ssec:k2_tri}
    Setting $k=2$ we can find that the modal correction mass matrix is:
    \begin{equation}
        \tb{Q} = \begin{bmatrix}
            0 & 0 & 0 & 0 & 0 & 0 \\
            0 & 0 & 0 & 0 & 0 & 0 \\
            0 & 0 & (5q_0 - q_1)/4 & 0 & 0 & (q_0 - q_1)\sqrt{5/4} \\
            0 & 0 & 0 & 0 & 0 & 0 \\
            0 & 0 & 0 & 0 & q_1 & 0 \\
            0 & 0 & (q_0 - q_1)\sqrt{5/4} & 0 & 0 & q_0
        \end{bmatrix}.
    \end{equation}
    From \cref{eq:la_cond6} we have the requirement that $\tb{M}+\tb{Q}$ is positive definite. This implies that the leading diagonal of a Cholesky factorisation of $\tb{M}+\tb{Q}$ has to be real and positive, leading to the following conditions on stability:
    \begin{equation}
            q_1 > -1 \quad \mathrm{and} \quad  9q_0 - 5q_1 > -4.
    \end{equation}
    The single parameter FR scheme of \cref{def:cast_fr} is then recovered for 
    \begin{equation}
            q_0 = \frac{410c}{3}, \quad \mathrm{and} \quad q_1 = 150c,
    \end{equation}
    which leads to the stability condition that
    \begin{equation}
        c > -\frac{1}{150}.
    \end{equation}
 
\subsection{$k=3$}
    Setting $k=3$ we can find that the modal correction mass matrix is:
    \begin{equation}
        \tb{Q} = \left[ \begin{array}{cccccccc}
            \mathbf{0}  & \multicolumn{7}{c}{\mathbf{0}} \\
      \multirow{7}{*}{$\boldsymbol{0}$} & (q_1 + 4q_0)/5 + q_2\sqrt{48/175} & 0 & 0 & 0 & 0 & (q_0 - q_1)/\sqrt{5} + q_2\sqrt{3/35} & 0 \\
      & 0 & 0 & 0 & 0 & 0 & 0 & 0 \\
      & 0 & 0 & 0 & 0 & 0 & 0 & 0 \\
      & 0 & 0 & 0 & q_0 - q_2\sqrt{16/21} & 0 & 0 & q_2 \\
      & 0 & 0 & 0 & 0 & 0 & 0 & 0 \\
      & (q_0 - q_1)/\sqrt{5} + q_2\sqrt{3/35} & 0 & 0 & 0 & 0 & q_1 & 0 \\
      & 0 & 0 & 0 & q_2 & 0 & 0 & q_0
        \end{array}\right]
    \end{equation}
    with stability limits of:
    \begin{subequations}
        \begin{align}
            28q_0 + 7q_1 + 4\sqrt{21}q_2 &> -35,\\
            21q_0 - 4\sqrt{21}q_2 &> -21,\\
            (7 q_0 + \sqrt{21}q_2 + 7)(7q_0 - 42q_1 + \sqrt{21} q_2 - 35) &> 0,\\
            (7q_0 + \sqrt{21}q_2 + 7)(3q_0 - \sqrt{21}q_2 + 3) &>0.
        \end{align}
    \end{subequations}
    The single parameter scheme of \cref{def:cast_fr} is then recovered with
    \begin{equation}
        q_0 = 6384c,\quad  q_1 = \frac{27440c}{3},\quad \mathrm{and} \quad q_2 = -168\sqrt{21}c,
    \end{equation}
    subject to the stability condition that
    \begin{equation}
        c > -\frac{1}{9800}.
    \end{equation}
    
\begin{figure}[tbhp]
        \centering
        \subfloat[Flux point 1.]{\label{fig:k3_tri_dg_f1}\adjustbox{width=0.38\linewidth,valign=b}{    \begin{tikzpicture}[scale=2]
		\begin{axis}[name=plot1,xlabel={$x$},ylabel={$y$},
    		axis line style={latex-latex},
            axis y line=middle,
            axis x line=middle,
            xmode=linear, 
            ymode=linear, 
		    xtick={0},
		    ytick={0},
		    xticklabels={ },
		    yticklabels={ },
		    ylabel style={rotate=0, shift={(-0.09,0.06)}},
		    xlabel style={rotate=0, shift={(0.03,0.01)}},
    		xmin=-1.1,xmax=1.1,
    		ymin=-0.6351,ymax=1.2702,
    		width=5cm,
    		height=4.330125cm,
			scale only axis=true,
    		style={font=\normalsize},
    		axis on top,
        ]
            \addplot[thick] graphics[xmin=-1,ymin=-0.5774,xmax=1,ymax=1.1547] {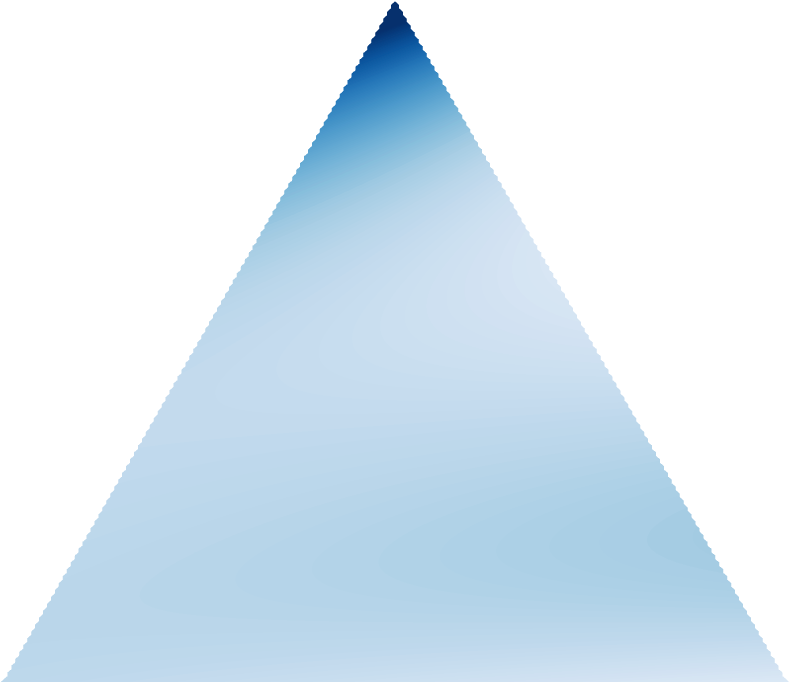};
            
            \draw[black] (axis cs:-1,-0.5774) -- (axis cs:1,-0.5774) -- (axis cs:0,1.1547) -- cycle;
            
     \draw[Set1-A,fill=Set1-A] (axis cs:-0.069432,1.0344) circle (0.018);
     \draw[black,fill=black] (axis cs:-0.33001,0.58311) circle (0.018);
     \draw[black,fill=black] (axis cs:-0.66999,-0.0057571) circle (0.018);
     \draw[black,fill=black] (axis cs:-0.93057,-0.45709) circle (0.018);
     \draw[black,fill=black] (axis cs:0.93057,-0.45709) circle (0.018);
     \draw[black,fill=black] (axis cs:0.66999,-0.0057571) circle (0.018);
     \draw[black,fill=black] (axis cs:0.33001,0.58311) circle (0.018);
     \draw[black,fill=black] (axis cs:0.069432,1.0344) circle (0.018);
     \draw[black,fill=black] (axis cs:-0.86114,-0.57735) circle (0.018);
     \draw[black,fill=black] (axis cs:-0.33998,-0.57735) circle (0.018);
     \draw[black,fill=black] (axis cs:0.33998,-0.57735) circle (0.018);
     \draw[black,fill=black] (axis cs:0.86114,-0.57735) circle (0.018);
         
		\end{axis}
    \begin{scope}[on behind layer]
        \draw[black!00] (0,-0.265) rectangle (5,4.4);
    \end{scope}
\end{tikzpicture}}}
        ~
        \subfloat[]{\label{fig:k3_tri_dg_f2}\adjustbox{width=0.55\linewidth,valign=b}{    \begin{tikzpicture}[scale=2]
		\begin{axis}[name=plot1,xlabel={$x$},ylabel={$y$},
    		axis line style={latex-latex},
            axis y line=middle,
            axis x line=middle,
            xmode=linear, 
            ymode=linear, 
		    xtick={0},
		    ytick={0},
		    xticklabels={ },
		    yticklabels={ },
		    ylabel style={rotate=0, shift={(-0.09,0.06)}},
		    xlabel style={rotate=0, shift={(0.03,0.01)}},
    		xmin=-1.1,xmax=1.1,
    		ymin=-0.6351,ymax=1.2702,
    		width=5cm,
    		height=4.330125cm,
			scale only axis=true,
    		style={font=\normalsize},
    		axis on top,
    		colorbar,
    		colormap/Blues,
    		colorbar style={
    		    point meta max=10, point meta min=-4,
                ylabel=$\nabla\cdot\mathbf{h}$,
                ylabel style={font=\normalsize},
                ytick={-4,-2,0,2,4,6,8,10},
                yticklabel style={
                    text width=2em,
                    align=left,
                    font=\normalsize,
                }
            }
        ]
            \addplot[thick] graphics[xmin=-1,ymin=-0.5774,xmax=1,ymax=1.1547] {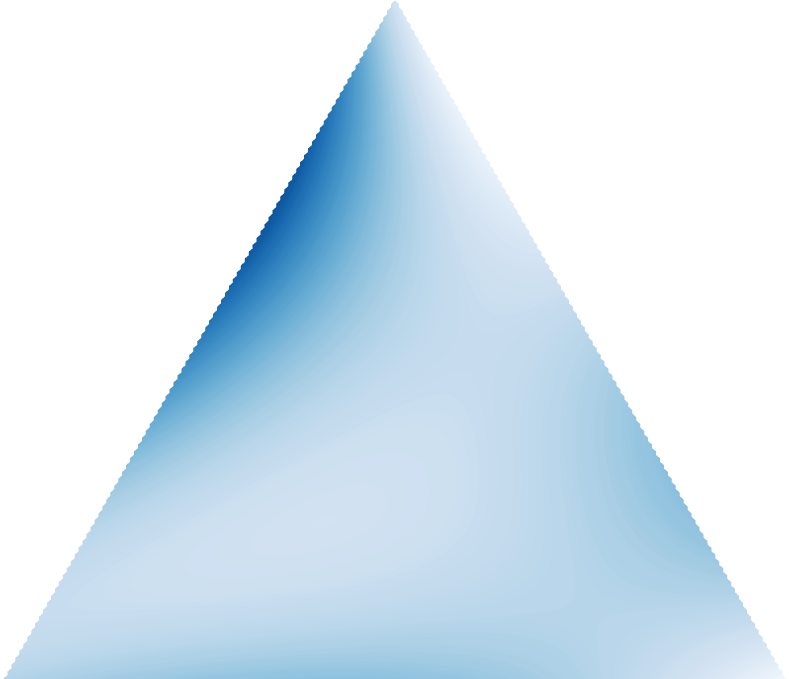};
            
            \draw[black] (axis cs:-1,-0.5774) -- (axis cs:1,-0.5774) -- (axis cs:0,1.1547) -- cycle;
            
     \draw[black,fill=black] (axis cs:-0.069432,1.0344) circle (0.018);
     \draw[Set1-A,fill=Set1-A] (axis cs:-0.33001,0.58311) circle (0.018);
     \draw[black,fill=black] (axis cs:-0.66999,-0.0057571) circle (0.018);
     \draw[black,fill=black] (axis cs:-0.93057,-0.45709) circle (0.018);
     \draw[black,fill=black] (axis cs:0.93057,-0.45709) circle (0.018);
     \draw[black,fill=black] (axis cs:0.66999,-0.0057571) circle (0.018);
     \draw[black,fill=black] (axis cs:0.33001,0.58311) circle (0.018);
     \draw[black,fill=black] (axis cs:0.069432,1.0344) circle (0.018);
     \draw[black,fill=black] (axis cs:-0.86114,-0.57735) circle (0.018);
     \draw[black,fill=black] (axis cs:-0.33998,-0.57735) circle (0.018);
     \draw[black,fill=black] (axis cs:0.33998,-0.57735) circle (0.018);
     \draw[black,fill=black] (axis cs:0.86114,-0.57735) circle (0.018);
         
		\end{axis}
\end{tikzpicture}}}
        \caption{\label{fig:k3_dg_corr}Divergence of the correction field for $k=3$ with $q_0=q_1=q_2=0$, \ie DG, for two flux points shown in red.}
    \end{figure}

\begin{figure}[tbhp]
        \centering
        \subfloat[Flux point 1.]{\label{fig:k3_tri_q_f1}\adjustbox{width=0.38\linewidth,valign=b}{    \begin{tikzpicture}[scale=2]
		\begin{axis}[name=plot1,xlabel={$x$},ylabel={$y$},
    		axis line style={latex-latex},
            axis y line=middle,
            axis x line=middle,
            xmode=linear, 
            ymode=linear, 
		    xtick={0},
		    ytick={0},
		    xticklabels={ },
		    yticklabels={ },
		    ylabel style={rotate=0, shift={(-0.09,0.06)}},
		    xlabel style={rotate=0, shift={(0.03,0.01)}},
    		xmin=-1.1,xmax=1.1,
    		ymin=-0.6351,ymax=1.2702,
    		width=5cm,
    		height=4.330125cm,
			scale only axis=true,
    		style={font=\normalsize},
    		axis on top,
        ]
            \addplot[thick] graphics[xmin=-1,ymin=-0.5774,xmax=1,ymax=1.1547] {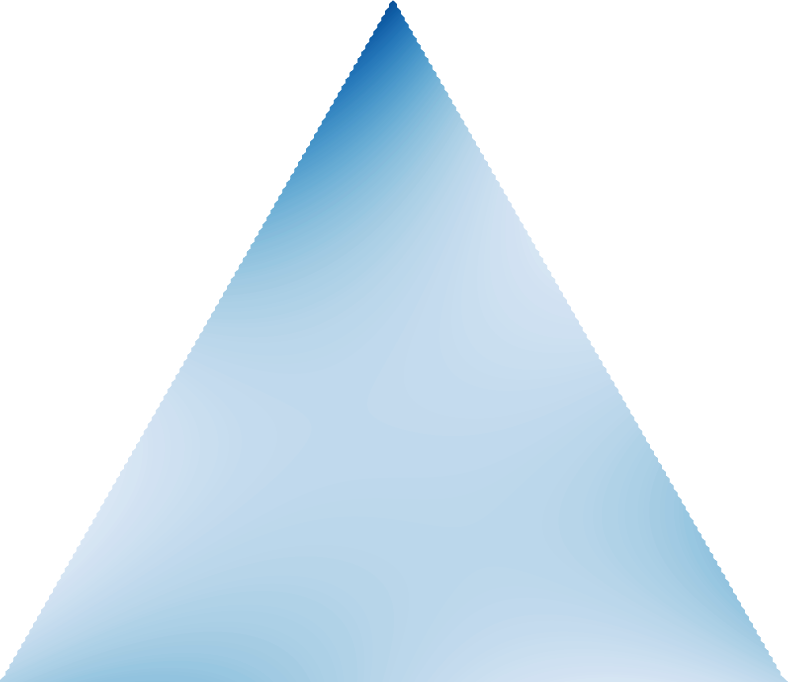};
            
            \draw[black] (axis cs:-1,-0.5774) -- (axis cs:1,-0.5774) -- (axis cs:0,1.1547) -- cycle;
            
     \draw[Set1-A,fill=Set1-A] (axis cs:-0.069432,1.0344) circle (0.018);
     \draw[black,fill=black] (axis cs:-0.33001,0.58311) circle (0.018);
     \draw[black,fill=black] (axis cs:-0.66999,-0.0057571) circle (0.018);
     \draw[black,fill=black] (axis cs:-0.93057,-0.45709) circle (0.018);
     \draw[black,fill=black] (axis cs:0.93057,-0.45709) circle (0.018);
     \draw[black,fill=black] (axis cs:0.66999,-0.0057571) circle (0.018);
     \draw[black,fill=black] (axis cs:0.33001,0.58311) circle (0.018);
     \draw[black,fill=black] (axis cs:0.069432,1.0344) circle (0.018);
     \draw[black,fill=black] (axis cs:-0.86114,-0.57735) circle (0.018);
     \draw[black,fill=black] (axis cs:-0.33998,-0.57735) circle (0.018);
     \draw[black,fill=black] (axis cs:0.33998,-0.57735) circle (0.018);
     \draw[black,fill=black] (axis cs:0.86114,-0.57735) circle (0.018);
         
		\end{axis}
    \begin{scope}[on behind layer]
        \draw[black!00] (0,-0.265) rectangle (5,4.4);
    \end{scope}
\end{tikzpicture}}}
        ~
        \subfloat[Flux point 2.]{\label{fig:k3_tri_q_f2}\adjustbox{width=0.55\linewidth,valign=b}{    \begin{tikzpicture}[scale=2]
		\begin{axis}[name=plot1,xlabel={$x$},ylabel={$y$},
    		axis line style={latex-latex},
            axis y line=middle,
            axis x line=middle,
            xmode=linear, 
            ymode=linear, 
		    xtick={0},
		    ytick={0},
		    xticklabels={ },
		    yticklabels={ },
		    ylabel style={rotate=0, shift={(-0.09,0.06)}},
		    xlabel style={rotate=0, shift={(0.03,0.01)}},
    		xmin=-1.1,xmax=1.1,
    		ymin=-0.6351,ymax=1.2702,
    		width=5cm,
    		height=4.330125cm,
			scale only axis=true,
    		style={font=\normalsize},
    		axis on top,
    		colorbar,
    		colormap/Blues,
    		colorbar style={
    		    point meta max=10, point meta min=-4,
                ylabel=$\nabla\cdot\mathbf{h}$,
                ylabel style={font=\normalsize},
                ytick={-4,-2,0,2,4,6,8,10},
                yticklabel style={
                    text width=2em,
                    align=left,
                    font=\normalsize,
                }
            }
        ]
            \addplot[thick] graphics[xmin=-1,ymin=-0.5774,xmax=1,ymax=1.1547] {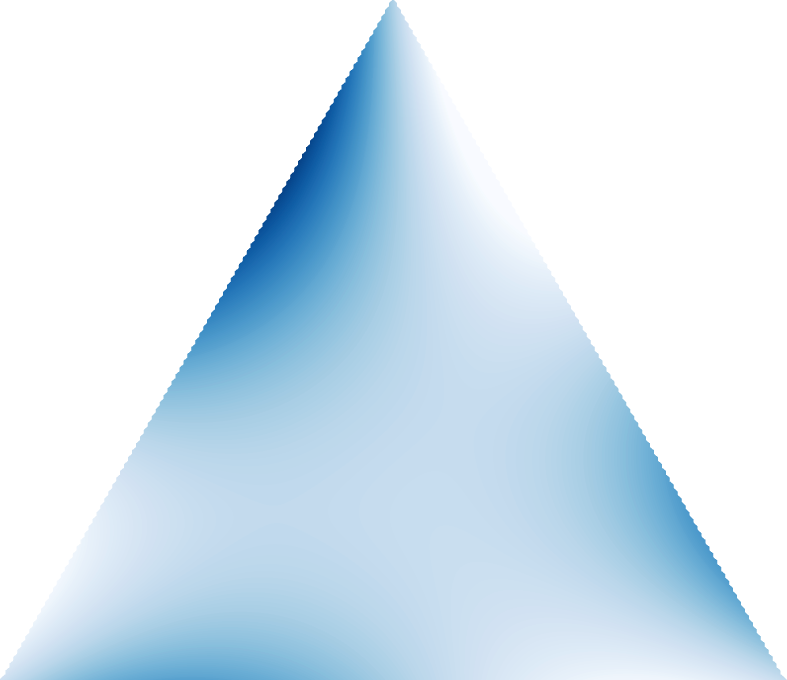};
            
            \draw[black] (axis cs:-1,-0.5774) -- (axis cs:1,-0.5774) -- (axis cs:0,1.1547) -- cycle;
            
            \draw[black,fill=black] (axis cs:-0.069432,1.0344) circle (0.018);
            \draw[Set1-A,fill=Set1-A] (axis cs:-0.33001,0.58311) circle (0.018);
            \draw[black,fill=black] (axis cs:-0.66999,-0.0057571) circle (0.018);
            \draw[black,fill=black] (axis cs:-0.93057,-0.45709) circle (0.018);
            \draw[black,fill=black] (axis cs:0.93057,-0.45709) circle (0.018);
            \draw[black,fill=black] (axis cs:0.66999,-0.0057571) circle (0.018);
            \draw[black,fill=black] (axis cs:0.33001,0.58311) circle (0.018);
            \draw[black,fill=black] (axis cs:0.069432,1.0344) circle (0.018);
            \draw[black,fill=black] (axis cs:-0.86114,-0.57735) circle (0.018);
            \draw[black,fill=black] (axis cs:-0.33998,-0.57735) circle (0.018);
            \draw[black,fill=black] (axis cs:0.33998,-0.57735) circle (0.018);
            \draw[black,fill=black] (axis cs:0.86114,-0.57735) circle (0.018);
         
		\end{axis}
\end{tikzpicture}}}
        \caption{\label{fig:k3_q_corr}Divergence of the correction field for $k=3$ with $q_0=q_1=q_2=1$ for two flux points shown in red.}
    \end{figure}

    As an example of how the correction functions are effected by $\tb{Q}$, consider \cref{fig:k3_dg_corr} which shows the divergence of the DG correction field and \cref{fig:k3_q_corr} which shows the divergence of the correction field for $\tb{Q}(q_0=1,q_1=1,q_2=1)$.

\subsection{$k=4$}
    Finally repeating this analysis for $k=4$, we find the following definition of the $\tb{Q}$ matrix.
    \begin{subequations}
        \begin{align}
            \tb{Q} &= \left[ \begin{array}{cccccccccccc}
                \mathbf{0} & \multicolumn{11}{c}{\mathbf{0}} \\
                \multirow{11}{*}{$\boldsymbol{0}$}
                  & \theta_4 & 0 & 0 & 0 & 0 & 0 & 0 & q_3 & 0 & 0 & q_2 \\
          & 0 & 0 & 0 & 0 & 0 & 0 & 0 & 0 & 0 & 0 & 0 \\
          & 0 & 0 & 0 & 0 & 0 & 0 & 0 & 0 & 0 & 0 & 0 \\
          & 0 & 0 & 0 & 0 & 0 & 0 & 0 & 0 & 0 & 0 & 0 \\
          & 0 & 0 & 0 & 0 & q_1 & 0 & 0 & 0 & 0 & q_0 & 0 \\
          & 0 & 0 & 0 & 0 & 0 & 0 & 0 & 0 & 0 & 0 & 0 \\
          & 0 & 0 & 0 & 0 & 0 & 0 & 0 & 0 & 0 & 0 & 0 \\
          & q_3 & 0 & 0 & 0 & 0 & 0 & 0 & \theta_3 & 0 & 0 & \theta_2 \\
          & 0 & 0 & 0 & 0 & 0 & 0 & 0 & 0 & 0 & 0 & 0 \\
          & 0 & 0 & 0 & 0 & q_0 & 0 & 0 & 0 & 0 & \theta_1 & 0 \\
          & q_2 & 0 & 0 & 0 & 0 & 0 & 0 & \theta_2 & 0 & 0 & \theta_0
            \end{array} \right]\\
            \mathrm{for} & \begin{cases}
                \theta_0 &= (16\sqrt{21}q_0 + 98q_1 + 99q_2 - 9\sqrt{5}q_3)/98 \\
                \theta_1 &= (8\sqrt{21}q_0 + 42q_1 + 75q_2 - 45\sqrt{5}q_3)/42 \\
                \theta_2 &= (5\sqrt{105}q_0 + 20\sqrt{5}q_2 - 18q_3)/49 \\
                \theta_3 &= (4\sqrt{21}q_0 + 147q_1 + 240q_2  - 102\sqrt{5}q_3)/147 \\
                \theta_4 &= (2\sqrt{21}q_0 + 6q_1 + 12q_2 - 3\sqrt{5}q_3)/6
            \end{cases}
        \end{align}
    \end{subequations}
    subject to the constraints on stability that:
    \begin{subequations}
        \begin{align}
            \theta_4 &> -1 \\
            q_1 &> -1\\
            - q_3^2 + \theta_3 + \theta_4 + \theta_3\theta_4 &> -1 \\
            - q_0^2 + q_1 + \theta_1 + q_1\theta_1 &> -1 \\
            2 \theta_2q_2 q_3 + \theta_3 (\theta_4 - q_2^2 + 1) + \theta_0 ((\theta_3 + 1) (\theta_4 + 1) - q_3^2) - (\theta_4 + 1) \theta_2^2 + \theta_4 - q_2^2 - q_3^2 &> -1.
        \end{align}
    \end{subequations}
    The single parameter family of \cref{def:cast_fr} is then found as a subset when
    \begin{equation}
        q_0 = -30240\sqrt{21}c,\quad q_1 = 56440c,\quad q_2 = 5040c, \quad \mathrm{and} \quad q_3 = -80640\sqrt{5}c,
    \end{equation}
    with the stability condition of:
    \begin{equation}
        c > \frac{\sqrt{1129}-115}{76204800}.
    \end{equation}
    
    The procedure to find $\tb{Q}$ and the stability conditions can be generalised for arbitrary order using a symbolic manipulation tool. Performing analysis for higher orders, the stability conditions for the single parameter \citet{Castonguay2011} can be tabulated as in \cref{tab:cast_cm}. For $k>5$ it is does not seem possible to get a closed expression for the stability limit as the value of $c$ is the root of a high order polynomial. For example, at $k=6$ the polynomial is order seven. 
    
    \begin{table}[tbhp]
        \centering
        \begin{tabular}{c c}
            \toprule
            $k$ & $c>$ \\ \midrule
            1 & $-\frac{1}{6}$\\[4pt]
            2 & $-\frac{1}{150}$\\[4pt]
            3 & $-\frac{1}{9800}$\\[4pt]
            4 & $-\frac{115 - \sqrt{1129}}{76204800}$ \\[4pt]
            5 & $-\frac{67 - \sqrt{889}}{5488560000}$\\[4pt]
            6 & $\num{-2.88363e-11}$\\[4pt]
           \bottomrule
        \end{tabular}
        \caption{\label{tab:cast_cm}Stability limits for single parameter \citet{Castonguay2011} correction functions on triangles.}
    \end{table}
    
\section{Spectral Difference Methods}\label{sec:sd}
    The spectral difference (SD) method~\citep{Liu2006} is a high-order method similar to FR but with the flux evaluated at staggered set of points, analogous to the method of \citet{Kopriva1996}. The nodal basis of the flux is then chosen such that it lies in the Raviart--Thomas~\citep{Raviart1977} space of the approximation space \footnote{For this reason SD is sometimes referred to as Raviart--Thomas SD.}. This has the effect of elevating the flux function order, which has been found to give rise to better aliasing properties~\citep{Cox2021}.
    
    In the 1D linear case, SD can be found to be a member of the one parameter class of FR methods~\citep{Vincent2010,Trojak2021}. Using the generalisation of \citet{Trojak2021}, the SD correction functions can be expressed as
    \begin{equation}
        h_L(x) = \frac{(-1)^k}{2}(1-x)J^{(\alpha,\beta)}_{k-1}(x), \quad \mathrm{and} \quad h_R(x) = \frac{1}{2}(x+1)J^{(\alpha,\beta)}_{k-1}(x).
    \end{equation}
    \citet{Jameson2010} has previously shown that in 1D the only linearly stable SD scheme is that corresponding to $(\alpha,\beta)=(0,0)$, \ie the interior flux points are located at the Gauss--Legendre nodes. \citet{Trojak2021} showed that Fourier analysis can be effectively used to find stable SD schemes with alternative point layouts for which linear stability proofs could not be constructed.
    
    A long-standing difficulty has existed in defining linearly stable SD schemes on triangles. Schemes can be constructed for tensor product elements on a maximal order basis, such as quadrilaterals and hexahedrals. However, simplex elements have proven to be more difficult, with some schemes found via Fourier analysis that are stable under interface upwinding. The broad set of stable schemes outlined in \cref{sec:erfr_tri} offers a promising route to find generally stable SD method. 

\subsection{One-dimension}\label{sec:sd_1d}
    As an initial test of a procedure to find SD correction functions, we look to confirm the SD stability theorem of \citet{Jameson2010} in 1D. Here we assume that the interior flux points are placed symmetrically within an element, and when there is an odd number of interior points a single point is placed at the centre. A numerical version of this study has previously been performed by \citet{Abeele2008}.
    
    \begin{figure}[tbhp]
        \centering
        \adjustbox{width=0.6\linewidth, valign=b}{\begin{tikzpicture}[scale=1]

    \draw[black!00] (-5.5,-1.5) rectangle (5.5,1.5);

    \begin{scope}[on behind layer]
        \draw[black, thick] (-5,0) -- (5,0);
        \draw[black, thick] ( 5,-0.25) -- ( 5,0.25);
        \draw[black, thick] (-5,-0.25) -- (-5,0.25);
    \end{scope}
    
    \begin{scope}
        \fill[Set1-C] (-5,0) circle (3pt);
        \fill[Set1-C] ( 5,0) circle (3pt);
    \end{scope}
    
    \begin{scope}
        \fill[Set1-A] (-3.5,0) circle (3pt);
        \fill[Set1-A] ( 3.5,0) circle (3pt);
        \draw[black, -latex, thick] (0,-0.75) -- (-3.5,-0.75) node[midway,below] {$-z_1$};
        \draw[black, |-latex, thick] (0,-0.75) -- ( 3.5,-0.75) node[midway,below] {$z_1$};
        
        \fill[Set1-B] (-1.5,0) circle (3pt);
        \fill[Set1-B] ( 1.5,0) circle (3pt);
        \draw[black, -latex, thick] (0,0.75) -- (-1.5,0.75) node[midway,above] {$-z_2$};
        \draw[black, |-latex, thick] (0,0.75) -- ( 1.5,0.75) node[midway,above] {$z_2$};
    \end{scope}
    
\end{tikzpicture}}
        \caption{\label{fig:sd_1d}}
    \end{figure}
    
    Using a point layout similar to the example shown in \cref{fig:sd_1d}, the nodal representation of the 1D correction functions can be written as 
    \begin{equation}
        h_L(x) = \frac{(1-x)x^{n}\prod^{m}_{i=1}(x-z_i)(x+z_i)}{(-1)^{p+1}2\prod^{m}_{i=1}(1+z_i)(1-z_i)}, \quad \mathrm{and} \quad
        h_R(x) = \frac{(1+x)x^{n}\prod^{m}_{i=1}(x-z_i)(x+z_i)}{2\prod^{m}_{i=1}(1+z_i)(1-z_i)},
    \end{equation}
    for $m=\floor{p/2}$ and $n=p\mod{2}$. This can then be differentiated and transformed into a modal representation allowing for $\tb{Q}$ to be found via:
    \begin{equation}
        \tb{Q}\tb{C}_\mathrm{SD} = -\tb{M}(\tb{C}_\mathrm{SD} - \tb{C}_\mathrm{DG}).
    \end{equation}
    Here the interpolation operators have been factored out by using the DG correction matrix to give system that is more straightforwardly solved. 
    For $k=3$, we have $m=1$ and $n=1$, and we find $\tb{Q}$ as
    \begin{equation}
        \tb{Q} = \begin{bmatrix}
            0 & 0 & 0 & 0\\
            0 & 0 & 0 & q_1 \\
            0 & 0 & q_2 & 0 \\
            0 & q_1 & 0 & q_0
        \end{bmatrix}, \;\;\mathrm{for}
        \begin{cases}
            q_0 &= -\frac{5(105z_1^4-161z_1^2+54)}{112},\\[7pt]
            q_1 &= \frac{3(3-5z_1^2)}{4}, \\[7pt]
            q_2 &= \frac{3-5z_1^2}{5}.
        \end{cases}
    \end{equation}
    However, from \cref{eq:la_cond5} we have the additional constraint that $q_1 = -5q_2/3$, which can only be satisfied if $q_1=q_2=0$. This occurs when $z_1=\pm\sqrt{3/5}$, which when combined with the flux point at $x=0$, gives the interior flux points as the nodes of the Gauss--Legendre quadrature. This confirms the result of \citet{Jameson2010} and when $z_1$ is substituted into $q_0$ we find $q_0=3/14$, as reported by \citet{Vincent2015}.
    
    Repeating this for $k=4$, we find that
    \begin{equation}
        \tb{Q} = \begin{bmatrix}
            0 & 0 & 0 & 0 & q_5\\
            0 & 0 & 0 & q_4 & 0\\
            0 & 0 & q_3 & 0 & q_2\\
            0 & q_4 & 0 & q_1 & 0\\
            q_5 & 0 & q_2 & 0 & q_0
        \end{bmatrix},\;\;
        \mathrm{for} \;\begin{cases}
            q_0 &= \frac{7(z_1^2 + z_2^2) - 10}{4}q_2 + \frac{63z_1^2z_2^2 - 63(z_1^2 +z_2^2) + 55}{36}, \\[7pt]
            q_1 &= \frac{525(z_1^2+z_2^2)z_1^2z_2^2 - 175(z_1^4 + z_2^4) - 105(z_1^2 + z_2^2) - 560z_1^2z_2^2 + 198}{336},\\[7pt]
            q_2 &= -\frac{525(z_1^2 + z_2^2)z_1^2z_2^2 - 175(z_1^4 + z_2^4) - 25(z_1^2+z_2^2) - 800z_1^2z_2^2 + 146}{240},\\[7pt]
            q_3 &= \frac{20q_2 - 35z_1^2z_2^2 + 21(z_1^2+z_2^2) -15}{35(z_1^2 + z_2^2) - 50)}, \\[7pt]
            q_4 &= \frac{15z_1^2z_2^2 - 5(z_1^2 + z_2^2) + 3}{12},\\[7pt]
            q_5 &=0.
        \end{cases}
    \end{equation}
    Again using the condition of \cref{eq:la_cond5}, we find that $q_i=0 \;\forall i\in\{1,\dots,5\}$ which is achieved when
    \begin{equation*}
        z_1^2 = \frac{3}{7} - \frac{2}{7}\sqrt{\frac{6}{5}}, \quad \mathrm{and} \quad z_2^2 = \frac{3}{7} + \frac{2}{7}\sqrt{\frac{6}{5}}.
    \end{equation*}
    or vice versa. This again corresponds to the Gauss--Legendre quadrature and gives $q_0=8/45$, corresponding to \citet{Vincent2015}.
    \begin{remark}
        The point symmetry imposed and the irrelevance of the ordering of the zeros is why there are multiple solutions that give a valid $\tb{Q}$. The correction functions recovered from each is the same. This symmetry can be further identified from the form of the $q$ terms and their lack odd powers of $z$.
    \end{remark}
    
\subsection{Triangular elements}\label{sec:sd_2d_tri}
    Next we extend this procedure to triangles.
    From the work of \citet{Balan2012}, we start by defining the Raviart--Thomas (RT) space in two-dimensions as
    \begin{equation}
        \mathbb{RT}_k = \mathbb{Q}_{k}\otimes\begin{bmatrix}1&0\\0&1\end{bmatrix} + \begin{bmatrix}
            x\\y
        \end{bmatrix}(\mathbb{Q}_{k}-\mathbb{Q}_{k-1}),
    \end{equation}
    which for $k=2$ gives
    \begin{equation}
        \mathrm{span}(\mathbb{RT}_k) = \left\{\begin{pmatrix}\phi_1 \\ 0 \end{pmatrix}, \begin{pmatrix}\phi_2 \\ 0 \end{pmatrix},\dots,\begin{pmatrix}0 \\ \phi_1 \end{pmatrix}, \begin{pmatrix}0 \\ \phi_2 \end{pmatrix}, \dots, \begin{pmatrix}x\phi_3 \\ y\phi_3 \end{pmatrix}, \begin{pmatrix}x\phi_5 \\ y\phi_5 \end{pmatrix},\begin{pmatrix}x\phi_6 \\ y\phi_6 \end{pmatrix} \right\}.
    \end{equation}
    In the FR method, the correction functions are within an RT space, and likewise the analogy of corrections in SD are within an RT space. For SD, this two-dimensional basis is then defined via a staggered or flux point set, $\{\pmb{\sigma}\}$, and requires normals to be associated with each point, $\mathbf{n}_s$. An example of these flux points and there normal can be seen in \cref{fig:sd_tri}. The Lagrange basis can then be defined via the Vandermonde as
    \begin{equation}
        \mathbf{V}_\mathrm{RT} = \begin{bmatrix} \mathbf{V}_{\mathrm{RT},x} & \mathbf{V}_{\mathrm{RT},y} \end{bmatrix}\cdot\mathbf{n}_s,
    \end{equation}
    where $\mathbf{V}_{\mathrm{RT},x}$ and $\mathbf{V}_{\mathrm{RT},x}$ are the Vandermonde matrices over $\mathbb{RT}_k\cdot[1,0]^T$ and $\mathbb{RT}_k\cdot[0,1]^T$ respectively. The Lagrange basis is then found from $\mathbf{V}_\mathrm{RT}^{-1}$. Finally, the corrections are set using this basis where, from the definition of the SD method, the trace of the SD flux points are located at the FR flux points.
    
\subsubsection{$k=1$}
    The most straightforward SD method to define on simplex element is for $k=1$. In this case a single interior flux point is required at the element centroid, with normals in $x$ and $y$. This case was not considered in \cref{sec:erfr_tri}, but the extended range of stable schemes can be found to be:
    \begin{equation}
        \tb{Q} = \begin{bmatrix}
            0 & 0 & 0\\
            0 & q_0 & 0\\
            0 & 0 & q_0
        \end{bmatrix}, \quad \mathrm{for} \quad q_0>-1.
    \end{equation}
    The SD correction function is then found to be recovered when $q_0=1/3$.
    
\subsubsection{$k=2$}
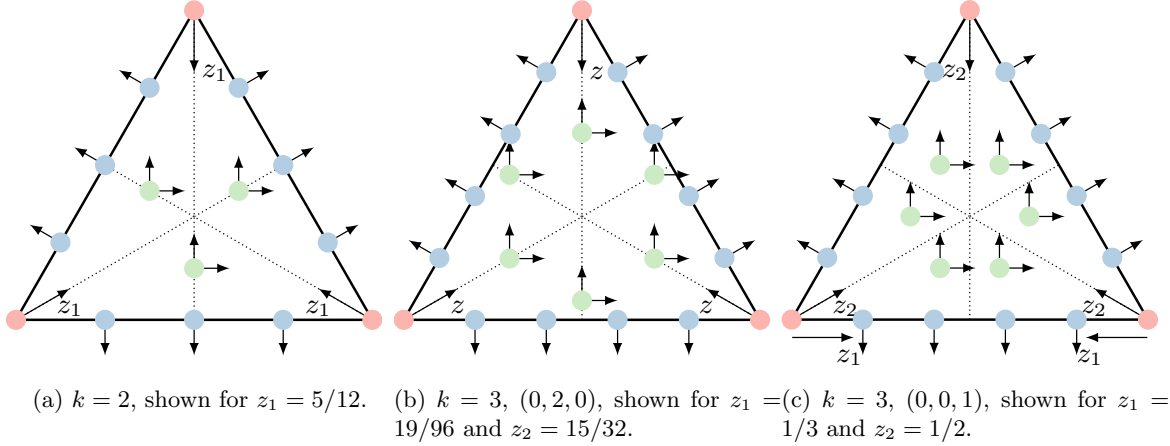
\begin{figure}[thbp]
        \centering
        \subfloat[$k=2$, shown for $z_1=5/12$.]{\label{fig:sd_tri_k2}\adjustbox{width=0.33\linewidth,valign=b}{\begin{tikzpicture}[scale=2]

    \begin{scope}
        \draw[fill={Pastel1-A}, color={Pastel1-A}] (-1,-0.57735026919) circle (1.5pt);
        
        \draw[fill={Pastel1-A}, color={Pastel1-A}] ( 1,-0.57735026919) circle (1.5pt);
        
        \draw[fill={Pastel1-A}, color={Pastel1-A}] ( 0, 1.15470053838) circle (1.5pt);
        
        \draw[fill={Pastel1-B},color={Pastel1-B}] (-0.75,-0.14434) circle (1.5pt);
        \draw[fill={Pastel1-B},color={Pastel1-B}] (-0.5,0.28868) circle (1.5pt);
        \draw[fill={Pastel1-B},color={Pastel1-B}] (-0.25,0.72169) circle (1.5pt);
        \draw[fill={Pastel1-B},color={Pastel1-B}] (0.25,0.72169) circle (1.5pt);
        \draw[fill={Pastel1-B},color={Pastel1-B}] (0.5,0.28868) circle (1.5pt);
        \draw[fill={Pastel1-B},color={Pastel1-B}] (0.75,-0.14434) circle (1.5pt);
        \draw[fill={Pastel1-B},color={Pastel1-B}] (0.5,-0.57735) circle (1.5pt);
        \draw[fill={Pastel1-B},color={Pastel1-B}] (0,-0.57735) circle (1.5pt);
        \draw[fill={Pastel1-B},color={Pastel1-B}] (-0.5,-0.57735) circle (1.5pt);

\draw[fill={Pastel1-C},color={Pastel1-C}] (-0.25,0.14434) circle (1.5pt);
\draw[fill={Pastel1-C},color={Pastel1-C}] (0,-0.28868) circle (1.5pt);
\draw[fill={Pastel1-C},color={Pastel1-C}] (0.25,0.14434) circle (1.5pt);
    \end{scope}
    
    \begin{scope}[on behind layer]
        \draw[thick, black] (-1,-0.57735026919) -- ( 1,-0.57735026919) -- ( 0, 1.15470053838) -- cycle;
                
        \draw[black,densely dotted] (-1,{-1/sqrt(3)}) -- (0.5,{0.5/sqrt(3)});
        \draw[black,densely dotted] (1,{-1/sqrt(3)}) -- (-0.5,{0.5/sqrt(3)});
        \draw[black,densely dotted] (0,{-1/sqrt(3)}) -- (0,{2/sqrt(3)});

        \draw[-latex,black] (-0.75,-0.14434) -- (-0.92321,-0.044338);
        \draw[-latex,black] (-0.5,0.28868) -- (-0.67321,0.38868);
        \draw[-latex,black] (-0.25,0.72169) -- (-0.42321,0.82169);
        \draw[-latex,black] (0.25,0.72169) -- (0.42321,0.82169);
        \draw[-latex,black] (0.5,0.28868) -- (0.67321,0.38868);
        \draw[-latex,black] (0.75,-0.14434) -- (0.92321,-0.044338);
        \draw[-latex,black] (0.5,-0.57735) -- (0.5,-0.77735);
        \draw[-latex,black] (0,-0.57735) -- (0,-0.77735);
        \draw[-latex,black] (-0.5,-0.57735) -- (-0.5,-0.77735);
        
        \draw[-latex,black] (-0.25,0.14434) -- (-0.05,0.14434);
        \draw[-latex,black] (0,-0.28868) -- (0.2,-0.28868);
        \draw[-latex,black] (0.25,0.14434) -- (0.45,0.14434);
        \draw[-latex,black] (-0.25,0.14434) -- (-0.25,0.34434);
        \draw[-latex,black] (0,-0.28868) -- (0,-0.088675);
        \draw[-latex,black] (0.25,0.14434) -- (0.25,0.34434);
        
        \draw[black,-latex] (1,{-1/sqrt(3)}) -- (0.6969,-0.4024) node[below] {\footnotesize{$z_1$}};
        \draw[black,-latex] (-1,{-1/sqrt(3)}) -- (-0.6969,-0.4024) node[below] {\footnotesize{$z_1$}};
        \draw[black,-latex] (0,{2/sqrt(3)}) -- (0,{2/sqrt(3) - 0.35}) node[right,shift={(-0.05,0)}] {\footnotesize{$z_1$}};
    \end{scope}
    
    \begin{scope}[on behind layer]
        \draw[black!00] (-1,-0.89) rectangle (1,1.1547);
    \end{scope}
    
\end{tikzpicture}}}
        \subfloat[$k=3$, $(0,2,0)$, shown for $z_1=19/96$ and $z_2=15/32$.]{\label{fig:sd_tri_k3_1}\adjustbox{width=0.33\linewidth,valign=b}{\begin{tikzpicture}[scale=2]

    \begin{scope}
        \draw[fill={Pastel1-A}, color={Pastel1-A}] (-1,-0.57735026919) circle (1.5pt);
        
        \draw[fill={Pastel1-A}, color={Pastel1-A}] ( 1,-0.57735026919) circle (1.5pt);
        
        \draw[fill={Pastel1-A}, color={Pastel1-A}] ( 0, 1.15470053838) circle (1.5pt);
    
        \draw[fill={Pastel1-B},color={Pastel1-B}] (-0.8,-0.23094) circle (1.5pt);
        \draw[fill={Pastel1-B},color={Pastel1-B}] (-0.6,0.11547) circle (1.5pt);
        \draw[fill={Pastel1-B},color={Pastel1-B}] (-0.4,0.46188) circle (1.5pt);
        \draw[fill={Pastel1-B},color={Pastel1-B}] (-0.2,0.80829) circle (1.5pt);
        \draw[fill={Pastel1-B},color={Pastel1-B}] (0.2,0.80829) circle (1.5pt);
        \draw[fill={Pastel1-B},color={Pastel1-B}] (0.4,0.46188) circle (1.5pt);
        \draw[fill={Pastel1-B},color={Pastel1-B}] (0.6,0.11547) circle (1.5pt);
        \draw[fill={Pastel1-B},color={Pastel1-B}] (0.8,-0.23094) circle (1.5pt);
        \draw[fill={Pastel1-B},color={Pastel1-B}] (0.6,-0.57735) circle (1.5pt);
        \draw[fill={Pastel1-B},color={Pastel1-B}] (0.2,-0.57735) circle (1.5pt);
        \draw[fill={Pastel1-B},color={Pastel1-B}] (-0.2,-0.57735) circle (1.5pt);
        \draw[fill={Pastel1-B},color={Pastel1-B}] (-0.6,-0.57735) circle (1.5pt);
        
        \draw[fill={Pastel1-C},color={Pastel1-C}] (0.40625,-0.23455) circle (1.5pt);
        \draw[fill={Pastel1-C},color={Pastel1-C}] (0,0.4691) circle (1.5pt);
        \draw[fill={Pastel1-C},color={Pastel1-C}] (-0.40625,-0.23455) circle (1.5pt);
        \draw[fill={Pastel1-C},color={Pastel1-C}] (-0.40625,0.23455) circle (1.5pt);
        \draw[fill={Pastel1-C},color={Pastel1-C}] (0,-0.4691) circle (1.5pt);
        \draw[fill={Pastel1-C},color={Pastel1-C}] (0.40625,0.23455) circle (1.5pt);
        
    \end{scope}
    
    \begin{scope}[on behind layer]
        \draw[thick, black] (-1,-0.57735026919) -- ( 1,-0.57735026919) -- ( 0, 1.15470053838) -- cycle;
    
        \draw[-latex,black] (0.40625,-0.23455) -- (0.60625,-0.23455);
        \draw[-latex,black] (0,0.4691) -- (0.2,0.4691);
        \draw[-latex,black] (-0.40625,-0.23455) -- (-0.20625,-0.23455);
        \draw[-latex,black] (-0.40625,0.23455) -- (-0.20625,0.23455);
        \draw[-latex,black] (0,-0.4691) -- (0.2,-0.4691);
        \draw[-latex,black] (0.40625,0.23455) -- (0.60625,0.23455);
        \draw[-latex,black] (0.40625,-0.23455) -- (0.40625,-0.034549);
        \draw[-latex,black] (0,0.4691) -- (0,0.6691);
        \draw[-latex,black] (-0.40625,-0.23455) -- (-0.40625,-0.034549);
        \draw[-latex,black] (-0.40625,0.23455) -- (-0.40625,0.43455);
        \draw[-latex,black] (0,-0.4691) -- (0,-0.2691);
        \draw[-latex,black] (0.40625,0.23455) -- (0.40625,0.43455);

        \draw[-latex,black] (-0.8,-0.23094) -- (-0.97321,-0.13094);
        \draw[-latex,black] (-0.6,0.11547) -- (-0.77321,0.21547);
        \draw[-latex,black] (-0.4,0.46188) -- (-0.57321,0.56188);
        \draw[-latex,black] (-0.2,0.80829) -- (-0.37321,0.90829);
        \draw[-latex,black] (0.2,0.80829) -- (0.37321,0.90829);
        \draw[-latex,black] (0.4,0.46188) -- (0.57321,0.56188);
        \draw[-latex,black] (0.6,0.11547) -- (0.77321,0.21547);
        \draw[-latex,black] (0.8,-0.23094) -- (0.97321,-0.13094);
        \draw[-latex,black] (0.6,-0.57735) -- (0.6,-0.77735);
        \draw[-latex,black] (0.2,-0.57735) -- (0.2,-0.77735);
        \draw[-latex,black] (-0.2,-0.57735) -- (-0.2,-0.77735);
        \draw[-latex,black] (-0.6,-0.57735) -- (-0.6,-0.77735);
                
        \draw[black,densely dotted] (-1,{-1/sqrt(3)}) -- (0.5,{0.5/sqrt(3)});
        \draw[black,densely dotted] (1,{-1/sqrt(3)}) -- (-0.5,{0.5/sqrt(3)});
        \draw[black,densely dotted] (0,{-1/sqrt(3)}) -- (0,{2/sqrt(3)});
        
        \draw[black,-latex] (1,{-1/sqrt(3)}) -- (0.6969,-0.4024) node[below] {\footnotesize{$z$}};
        \draw[black,-latex] (-1,{-1/sqrt(3)}) -- (-0.6969,-0.4024) node[below] {\footnotesize{$z$}};
        \draw[black,-latex] (0,{2/sqrt(3)}) -- (0,{2/sqrt(3) - 0.35}) node[right,shift={(-0.05,0)}] {\footnotesize{$z$}};
    \end{scope}

    \begin{scope}[on behind layer]
        \draw[black!00] (-1,-0.89) rectangle (1,1.1547);
    \end{scope}
    
\end{tikzpicture}}}
        \subfloat[$k=3$, $(0,0,1)$, shown for $z_1=1/3$ and $z_2=1/2$.]{\label{fig:sd_tri_k3_2}\adjustbox{width=0.33\linewidth,valign=b}{\begin{tikzpicture}[scale=2]

    \begin{scope}
        \draw[fill={Pastel1-A}, color={Pastel1-A}] (-1,-0.57735026919) circle (1.5pt);
        
        \draw[fill={Pastel1-A}, color={Pastel1-A}] ( 1,-0.57735026919) circle (1.5pt);
        
        \draw[fill={Pastel1-A}, color={Pastel1-A}] ( 0, 1.15470053838) circle (1.5pt);
    
        \draw[fill={Pastel1-B},color={Pastel1-B}] (-0.8,-0.23094) circle (1.5pt);
        \draw[fill={Pastel1-B},color={Pastel1-B}] (-0.6,0.11547) circle (1.5pt);
        \draw[fill={Pastel1-B},color={Pastel1-B}] (-0.4,0.46188) circle (1.5pt);
        \draw[fill={Pastel1-B},color={Pastel1-B}] (-0.2,0.80829) circle (1.5pt);
        \draw[fill={Pastel1-B},color={Pastel1-B}] (0.2,0.80829) circle (1.5pt);
        \draw[fill={Pastel1-B},color={Pastel1-B}] (0.4,0.46188) circle (1.5pt);
        \draw[fill={Pastel1-B},color={Pastel1-B}] (0.6,0.11547) circle (1.5pt);
        \draw[fill={Pastel1-B},color={Pastel1-B}] (0.8,-0.23094) circle (1.5pt);
        \draw[fill={Pastel1-B},color={Pastel1-B}] (0.6,-0.57735) circle (1.5pt);
        \draw[fill={Pastel1-B},color={Pastel1-B}] (0.2,-0.57735) circle (1.5pt);
        \draw[fill={Pastel1-B},color={Pastel1-B}] (-0.2,-0.57735) circle (1.5pt);
        \draw[fill={Pastel1-B},color={Pastel1-B}] (-0.6,-0.57735) circle (1.5pt);

        \draw[fill={Pastel1-C},color={Pastel1-C}] (-0.16667,0.28868) circle (1.5pt);
        \draw[fill={Pastel1-C},color={Pastel1-C}] (0.16667,-0.28868) circle (1.5pt);
        \draw[fill={Pastel1-C},color={Pastel1-C}] (-0.33333,0) circle (1.5pt);
        \draw[fill={Pastel1-C},color={Pastel1-C}] (-0.16667,-0.28868) circle (1.5pt);
        \draw[fill={Pastel1-C},color={Pastel1-C}] (0.33333,0) circle (1.5pt);
        \draw[fill={Pastel1-C},color={Pastel1-C}] (0.16667,0.28868) circle (1.5pt);

    \end{scope}
    
    \begin{scope}[on behind layer]
        \draw[thick, black] (-1,-0.57735026919) -- ( 1,-0.57735026919) -- ( 0, 1.15470053838) -- cycle;
    
        \draw[-latex,black] (-0.16667,0.28868) -- (0.033333,0.28868);
        \draw[-latex,black] (0.16667,-0.28868) -- (0.36667,-0.28868);
        \draw[-latex,black] (-0.33333,0) -- (-0.13333,0);
        \draw[-latex,black] (-0.16667,-0.28868) -- (0.033333,-0.28868);
        \draw[-latex,black] (0.33333,0) -- (0.53333,0);
        \draw[-latex,black] (0.16667,0.28868) -- (0.36667,0.28868);
        \draw[-latex,black] (-0.16667,0.28868) -- (-0.16667,0.48868);
        \draw[-latex,black] (0.16667,-0.28868) -- (0.16667,-0.088675);
        \draw[-latex,black] (-0.33333,0) -- (-0.33333,0.2);
        \draw[-latex,black] (-0.16667,-0.28868) -- (-0.16667,-0.088675);
        \draw[-latex,black] (0.33333,0) -- (0.33333,0.2);
        \draw[-latex,black] (0.16667,0.28868) -- (0.16667,0.48868);

        \draw[-latex,black] (-0.8,-0.23094) -- (-0.97321,-0.13094);
        \draw[-latex,black] (-0.6,0.11547) -- (-0.77321,0.21547);
        \draw[-latex,black] (-0.4,0.46188) -- (-0.57321,0.56188);
        \draw[-latex,black] (-0.2,0.80829) -- (-0.37321,0.90829);
        \draw[-latex,black] (0.2,0.80829) -- (0.37321,0.90829);
        \draw[-latex,black] (0.4,0.46188) -- (0.57321,0.56188);
        \draw[-latex,black] (0.6,0.11547) -- (0.77321,0.21547);
        \draw[-latex,black] (0.8,-0.23094) -- (0.97321,-0.13094);
        \draw[-latex,black] (0.6,-0.57735) -- (0.6,-0.77735);
        \draw[-latex,black] (0.2,-0.57735) -- (0.2,-0.77735);
        \draw[-latex,black] (-0.2,-0.57735) -- (-0.2,-0.77735);
        \draw[-latex,black] (-0.6,-0.57735) -- (-0.6,-0.77735);
                
        \draw[black,densely dotted] (-1,{-1/sqrt(3)}) -- (0.5,{0.5/sqrt(3)});
        \draw[black,densely dotted] (1,{-1/sqrt(3)}) -- (-0.5,{0.5/sqrt(3)});
        \draw[black,densely dotted] (0,{-1/sqrt(3)}) -- (0,{2/sqrt(3)});
        
        \draw[black,-latex] (1,{-1/sqrt(3)}) -- (0.6969,-0.4024) node[below] {\footnotesize{$z_2$}};
        \draw[black,-latex] (-1,{-1/sqrt(3)}) -- (-0.6969,-0.4024) node[below] {\footnotesize{$z_2$}};
        \draw[black,-latex] (0,{2/sqrt(3)}) -- (0,{2/sqrt(3) - 0.35}) node[right,shift={(-0.45,0)}] {\footnotesize{$z_2$}};
        
        \draw[black,-latex] (-1,{-1/sqrt(3)-0.1}) -- (-0.65,{-1/sqrt(3)-0.1}) node[below,shift={(-0.05,0)}] {\footnotesize{$z_1$}};
        \draw[black,-latex] (1,{-1/sqrt(3)-0.1}) -- (0.65,{-1/sqrt(3)-0.1}) node[below,shift={(0.05,0)}] {\footnotesize{$z_1$}};
    \end{scope}

    \begin{scope}[on behind layer]
        \draw[black!00] (-1,-0.89) rectangle (1,1.1547);
    \end{scope}
    
\end{tikzpicture}}}
        \caption{\label{fig:sd_tri}Interior and boundary flux point locations and normals for SD on triangles.}
    \end{figure}

    We now consider $k=2$, here the number of flux points is 15, but this can be reduced to 12 degrees of freedom by repeating the interior flux point with orthogonal normals~\citep{Balan2012}. Similar to the method used in one-dimension, we parameterise the interior flux point locations by $z_1$, which can be placed using Barycentric coordinates in a manner ensuring rotational symmetry, see \cref{fig:sd_tri_k2}.
    
    Using this construction, the following matrix can be formed:
    \begin{equation}
        \tb{A} = \tb{Q}\tb{C}_\mathrm{SD} + \tb{M}(\tb{C}_\mathrm{SD} - \tb{C}_\mathrm{DG}),
    \end{equation}
    where a $\tb{Q}$ compatible with \cref{ssec:k2_tri} is sought such that $\tb{A}=0$. It was shown by \citet{Balan2012} that the stability of the method is independent of the boundary flux point locations, at least for linear equations, and so to reduce the complexity of the resulting matrices we place these points in an equispaced configuration. Focusing on the value of $\tb{A}_{0,0}$ we find that
    \begin{equation}
        \tb{A}_{0,0} = \frac{2-3 \sqrt{2}}{18 \sqrt[4]{3} z_1 (2 z_1-1)}.
    \end{equation}
    It is clear that there is no value of $z_1\in[0,1/2]$ that can satisfy $\tb{A}_{0,0}=0$.
    
    The assumption of collocated interior flux points can be relaxed and a second parameter $z_2$ can be introduced. Repeating the procedure above now with two variables, likewise no pair of variables, $(z_1,z_2)$, can be found for a norm in this class for which the energy monotonically decays in time. For brevity, a full display of the contradictions encountered is not given.
    
\subsection{$k=3$}
    For $k=3$, assuming collocated interior flux points, there are six degrees of freedom. For symmetric placement of these points there are two possible choices of orbits $(0,2,0)$ and $(0,0,1)$, based on the work of \citet{Witherden2015}, \ie two three-point orbits (parameterised by $z_1$ and $z_2$) or one six-point orbit (also parameterised by $z_1$ and $z_2$). Examples of these orbits are shown in \cref{fig:sd_tri_k3_1,fig:sd_tri_k3_2}.
    
    Starting with the $(0,2,0)$ configuration, we again use the result of \citet{Balan2011} that stability is independent of the boundary flux point location and use equispaced boundary flux points. Forming the Raviart--Thomas space and finding $\tb{A}$, we find from the second column of $\tb{A}$ that 
    \begin{equation}
        z_1 = \frac{6 + \sqrt{15}}{21}, \quad z_2 = \frac{6 - \sqrt{15}}{21}, \quad q_0 = 3/5, \quad \mathrm{and} \quad q_2 = 0.
    \end{equation}
    and a second solution with $z_1$ and $z_2$ swapped. Substituting these into $\tb{A}$ and studying the first column of $\tb{A}$, we see this leads to the contradiction of 
    \begin{equation}
        45q_1 - 7 = 0 \quad \mathrm{and} \quad 15q_1 + 11 = 0.
    \end{equation}
    Therefore, there is no $k=3$ SD scheme with the interior flux points in the configuration $(0,2,0)$ that is a form of filtered DG.

    Repeating this for interior flux points in the orbit $(0,0,1)$, we find in the first column of $\tb{A}$ that 
    \begin{equation}
        \tb{A}_{9,0} = \sqrt[4]{\frac{49}{3}} \quad \mathrm{and} \quad \tb{A}_{6,0} = \sqrt[4]{3}.
    \end{equation}
    Clearly this does not satisfy the condition that $\tb{A}=0$, from which we can draw the conclusion that there is no stable $k=3$ SD scheme, in either $(0,2,0)$ or $(0,0,1)$, that is a form of filtered DG.
    
    \begin{remark}
        This instability is similar to a finding presented by \citet{Abeele2008}, however, in that work only Fourier analysis was used to explore stability. Stable schemes where found by \citet{Veilleux2022} and \citet{Balan2012} using Fourier analysis, however, in that analysis interface upwinding was required to find stable schemes. Therefore, they are not strictly linearly stable. 
    \end{remark}

    Summarising, we found a linearly stable SD scheme for $k=1$, but show that for $k=2$ and $k=3$ none exist in this set of stable FR methods. It is unlikely that at yet higher orders stable SD schemes will be found in this set of FR methods, and taken with previous results, such as those of \citet{Balan2012}, \citet{Abeele2008}, and \citet{Veilleux2022}, it is unlikely that linearly stable SD methods on triangles can be found at all without upwind stabilisation.
\section{Numerical Experiments}\label{sec:numeical}
\begin{figure}[tbhp]
        \centering
        \adjustbox{width=0.48\linewidth,valign=b}{    \begin{tikzpicture}[scale=2]
		\begin{axis}[name=plot1,xlabel={$x$},ylabel={$y$},
    		axis line style={latex-latex},
            axis y line=middle,
            axis x line=middle,
            xmode=linear, 
            ymode=linear, 
		    xtick={0},
		    ytick={0},
		    xticklabels={ },
		    yticklabels={ },
		    ylabel style={rotate=0, shift={(-0.09,0.06)}},
		    xlabel style={rotate=0, shift={(0.03,0.01)}},
    		xmin=-11.0,xmax=11.0,
    		ymin=-11.0,ymax=11.0,
    		width=7cm,
    		height=7cm,
			scale only axis=true,
    		style={font=\normalsize},
    		axis on top,
    		colorbar,
    		colormap={lingrn}{rgb255=(14,28,28) rgb255=(29,102,71)  rgb255=(21,150,21) rgb255=(142,204,61) rgb255=(255,251,230)},
    		colorbar style={
    		    point meta max=1, point meta min=0.5,
                ylabel=$\rho$,
                ylabel style={font=\normalsize},
                ytick={0.5,0.6,0.7,0.8,0.9,1},
                yticklabel style={
                    text width=2em,
                    align=left,
                    font=\normalsize,
                }
            }
        ]
            \addplot[thick] graphics[xmin=-10,ymin=-10,xmax=10,ymax=10] {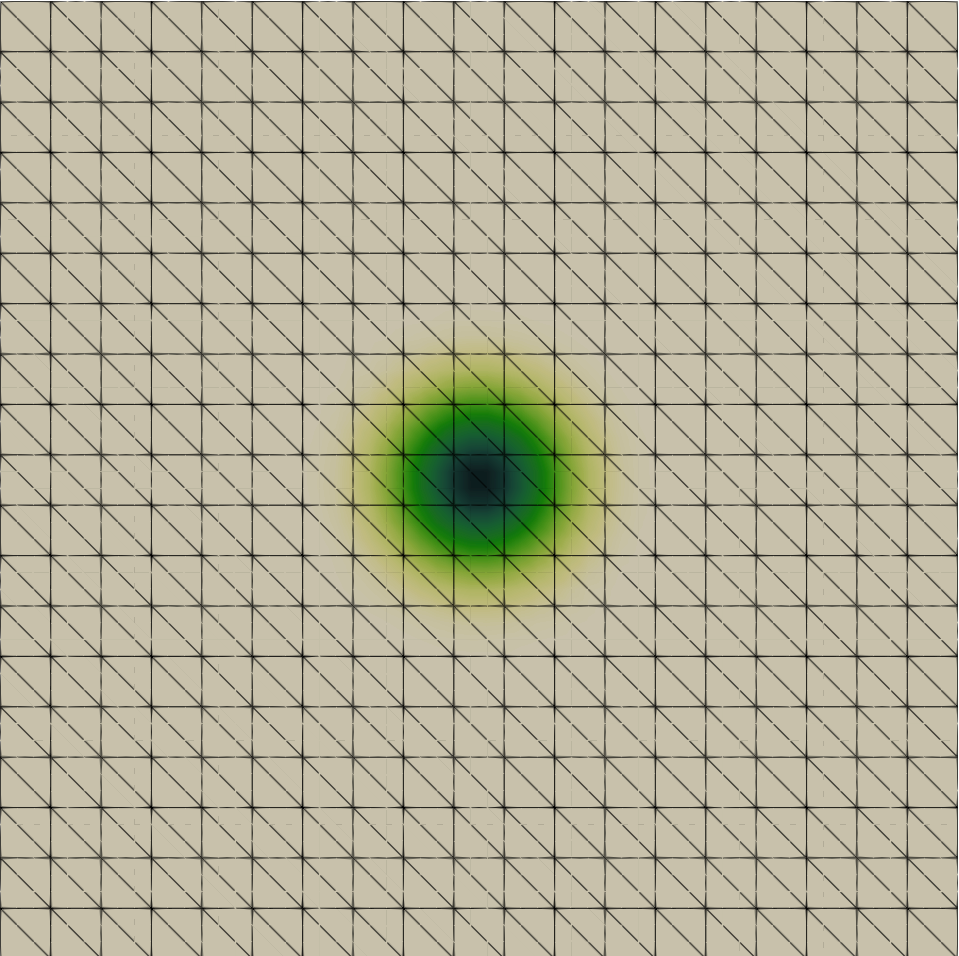};

		\end{axis}
\end{tikzpicture}}
        \caption{\label{fig:icv}Density contour for the Euler vortex with mesh shown for $n_x=20$.}
    \end{figure}

    To perform a numerical evaluation of the schemes defined here we considered the Euler vortex case~\citep{Shu1998}, a two-dimensional test case for the Euler equations. A periodic domain $\Omega=[-10,-10]^2$ subdivided into $2(n_x-1)^2$ regular right-angles triangles was used, see \cref{fig:icv}. The system of equations was then
    \begin{equation}
        \px{\mathbf{u}}{t} + \nabla\cdot\mathbf{F} = 0,\quad \mathrm{for} \quad
        \mathbf{u} = \begin{bmatrix}
            \rho \\ \boldsymbol{\rho v} \\ E
        \end{bmatrix}, \quad \mathrm{and} \quad \mathbf{F} = \begin{bmatrix}
            \boldsymbol{\rho v}\\
            \boldsymbol{\rho v}\otimes\mathbf{v} + P\mathbf{I}\\
            (E + P)\mathbf{v}
        \end{bmatrix},
    \end{equation}
    for pressure $P$ and energy $E$, and the initial condition was set as:
    \begin{subequations}
        \begin{align}
            \rho &= \left(1 - \half\left(\frac{\beta M}{\pi}\right)^2(\gamma - 1)\exp{(2\overline{r})}\right)^{1/(\gamma-1)},\\
            u &= \frac{\beta y\exp{(\overline{r})}}{2\pi R}, \\
            v &= 1 - \frac{\beta x\exp{(\overline{r})}}{2\pi R},\\
            P &= \frac{1}{\gamma M^2}\left(1 - \half\left(\frac{\beta M}{\pi}\right)^2(\gamma - 1)\exp{(2\overline{r})}\right)^{\gamma/(\gamma-1)},\\  
            \overline{r} &= \frac{1 - x^2 - y^2}{2R^2},
        \end{align}
    \end{subequations}
    where $M$ is the Mach number, $\beta$ is the vortex strength, and $R$ is the vortex width, set as $0.4$, $13.5$, and $1.5$ respectively. The error with time can then be calculated for a series of meshes, specifically we used the definition of $L_1$ and $L_2$ error of
    \begin{equation}
        E_1(t) = \int_K |\rho_\mathrm{exact}(t) - \rho(t)|\mathrm{d}\mathbf{x}, \quad \mathrm{and} \quad E_2(t) = \sqrt{\int_K (\rho_\mathrm{exact}(t) - \rho(t))^2\mathrm{d}\mathbf{x}},
    \end{equation}
    where the integrals are approximated with a degree 23 quadrature.
    
    In these tests, solution points were positioned at the quadrature points defined by \citet{Williams2014}. The common interface flux was calculated using a Rusanov flux and Einfeldt wavespeed predictions~\citep{Einfeldt1988} at flux points located with the Gauss--Legendre quadrature. For time integration, a standard explicit RK4 method was used.
    Results for $k=3$ are presented in \cref{tab:k3_tri_test} for $\mathbf{Q}=0$, $\mathbf{Q}_1(q_0=0.1,q_1=0.1,q_2=0.01)$, and $\mathbf{Q}_2(q_0=0,q_1=0,q_2=0.1)$. Here a constant time step of $\Delta t=\num{5e-3}$ was used and the $L_1$ and $L_2$ error is calculated at $t=100$. \cref{tab:k4_tri_test} shows the results for the test repeated for $k=4$, with $\mathbf{Q}=0$, $\mathbf{Q}_1(q_0=0.01,q_1=0.01,q_2=0.01,q_3=0.01)$, and $\mathbf{Q}_2(q_0=0.1,q_1=0,q_2=0,q_3=0)$. At $k=4$ a constant time step of $\Delta t=\num{2e-3}$ was used, again with error calculated at $t=100$. These data show that the correction functions tested were stable for $t\in[0,100]$ and furthermore the expected order of accuracy was recovered. The variation in error is evidence of the changing numerical properties caused by varying the correction function.
    
    \begin{table}[tbhp]
        \centering
        \begin{tabular}{c c c c c c c}
             \toprule
             \multirow{2}{*}{$n_x$} & \multicolumn{2}{c}{DG} & \multicolumn{2}{c}{$\mathbf{Q}_1$} & \multicolumn{2}{c}{$\mathbf{Q}_2$}\\\cmidrule(l){2-3} \cmidrule(l){4-5} \cmidrule(l){6-7} 
             & $E_1(t=100)$ & $E_2(t=100)$ & $E_1(t=100)$ & $E_2(t=100)$ & $E_1(t=100)$ & $E_2(t=100)$\\ 
             \midrule
             20 &  \num{5.57e-003}& \num{7.01e-003} & \num{5.77e-003} & \num{7.26e-003} & \num{5.65e-03} & \num{7.29e-03}\\
             25 &  \num{2.16e-003}& \num{2.65e-003} & \num{2.27e-003} & \num{2.79e-003} & \num{2.19e-03} & \num{2.73e-03}\\
             30 &  \num{1.07e-003}& \num{1.29e-003} & \num{1.11e-003} & \num{1.34e-003} & \num{1.10e-03} & \num{1.36e-03}\\
             35 &  \num{5.70e-004}& \num{6.61e-004} & \num{5.95e-004} & \num{6.96e-004} & \num{5.85e-04} & \num{6.96e-04}\\
             40 &  \num{3.23e-004}& \num{3.57e-004} & \num{3.36e-004} & \num{3.81e-004} & \num{3.36e-04} & \num{3.82e-04}\\ \midrule
             Order & 4.078 & 4.258 & 4.077 & 4.223 & 4.047 & 4.212\\\bottomrule
        \end{tabular}
        \caption{\label{tab:k3_tri_test}\noindent Error and order of the Euler vortex for $k=3$ FR with DG at $t=100$, $\mathbf{Q}_1(q_0=0.1,q_1=0.1,q_2=0.01)$, and $\mathbf{Q}_2(q_0=0,q_1=0,q_2=0.1)$.}
    \end{table}
    
    \begin{table}[tbhp]
        \centering
        \begin{tabular}{c c c c c c c}
             \toprule
             \multirow{2}{*}{$n_x$} & \multicolumn{2}{c}{DG} & \multicolumn{2}{c}{$\mathbf{Q}_1$} & \multicolumn{2}{c}{$\mathbf{Q}_2$}\\\cmidrule(l){2-3} \cmidrule(l){4-5} \cmidrule(l){6-7} 
             & $E_1(t=100)$ & $E_2(t=100)$ & $E_1(t=100)$ & $E_2(t=100)$ & $E_1(t=100)$ & $E_2(t=100)$\\\midrule
             20 &  \num{7.58e-04} & \num{8.39e-04} & \num{7.57e-04} & \num{8.46e-04} & \num{8.16e-04} & \num{9.75e-04}\\
             25 &  \num{2.40e-04} & \num{2.69e-04} & \num{2.41e-04} & \num{2.68e-04} & \num{2.66e-04}& \num{3.19e-04}\\
             30 &  \num{9.70e-05} & \num{1.10e-04} & \num{9.69e-05} & \num{1.11e-04} & \num{1.09e-04}& \num{1.34e-04}\\
             35 &  \num{4.50e-05} & \num{5.17e-05} & \num{4.49e-05} & \num{5.20e-05} & \num{5.11e-05}& \num{6.46e-05}\\
             40 &  \num{2.31e-05} & \num{2.70e-05} & \num{2.30e-05} & \num{2.71e-05} & \num{2.67e-05}& \num{3.40e-05}\\ \midrule
             Order & 5.028 & 4.951 & 5.033 & 4.952 & 4.931 & 4.830 \\\bottomrule
        \end{tabular}
        \caption{\label{tab:k4_tri_test}\noindent Error and order of the Euler vortex with $k=4$ FR for DG at $t=100$, $\mathbf{Q}_1(q_0=0.01,q_1=0.01,q_2=0.01,q_3=0.01)$, and $\mathbf{Q}_2(q_0=0.1,q_1=0,q_2=0,q_3=0)$.}
    \end{table}

\section{Conclusions}\label{sec:conclusions}
    A new multi-parameter set of stable flux reconstruction (FR) methods on triangles was constructed by using the summation-by-parts framework. The correction functions of \citet{Castonguay2011} were found to be a subset of this new stable set of FR methods, moreover we were able to successfully expand the stability region of \citet{Castonguay2011}. Using this new set of FR methods, we investigated if stable SD methods could be defined within it. We found that a stable SD scheme could be produced for $k=1$ and that none can be produced in this set of FR methods for $k=2$ and $k=3$. Numerical experiments were performed for a number of the correction functions outlined in this work and it was shown that the desired order of accuracy was recovered. The approaches outlined here can be used to find similar sets of methods on other element topologies which will be the subject of future work.


\section*{Data Availability}
The data that support the findings of this study are available from the corresponding author upon reasonable request.

\section*{Declaration}
\subsection*{Funding}
This work was supported by the Engineering and Physical Sciences Research Council (Grant number EP/R030340/1).

\subsection*{Competing Interests}
The authors have no relevant financial or non-financial interests to disclose.

\bibliographystyle{plainnat}
\bibliography{reference}


\begin{appendices}
\section{Weak Quadratures}\label{app:quad}
    To build the correction matrices defined in \cref{lem:lin_stab_q} a mass matrix is required, and it is often more practicable to produce this via quadrature rather than explicitly integrating the Lagrange basis. Yet, for triangular elements it is rarely possible to find a quadrature with both: sufficient strength to integrate the basis adequately and the same number of points as there are basis functions. For example, the quadrature rules of \citet{Williams2014} have $(k+1)(k+2)/2$ points but are not sufficiently accurate to integrate $\mathbb{Q}_{2k-1}$, as required for SBP to be valid. Instead, a more accurate quadrature like that of \citet{Witherden2015} can be used; however, these have more points. This poses a problem when calculating the correction function matrices within an implementation of FR as either the lumped mass matrix is insufficiently accurate or the Vandermonde matrix is not square.
    
    The solution is to use an $L_2$ projection with an intermediate set of points whose quadrature is sufficiently strong to integrate the basis. First consider the definition of the $L_2$ projection operators.
    \begin{definition}[$L_2$ projection]\label{def:l2_proj}
        For a nodal point set $\{\pmb{\zeta}_i\}_{i\leq N_q}$ defining some polynomial basis $Q\in\mathbb{Q}_q$ with an associated quadrature $\{\pmb{\omega}_i\}_{i\leq N_q}$ of strength at least $2q-1$, and a point set $\{\mathbf{x}_i\}_{i\leq N_k}$ defining $R\in\mathbb{Q}_k$ such that $N_q > N_k$ and $q>k$, then the $L_2$ projection matrix from $Q$ to $R$ is then
        \begin{equation}
            \mathbf{R}_{qk} = (\mathbf{P}_{qk}^T\mathbf{M}_q\mathbf{P}_{qk})^{-1}\mathbf{P}_{qk}^T\mathbf{M}_q,
        \end{equation}
        where $\mathbf{P}_{qk}$ is the prolongation matrix that interpolates from $R$ to $Q$.
    \end{definition}
    
    It is often more practical to set $\tb{Q}$ due to its sparser form. Therefore, we have the following lemma on the use of weaker quadratures:
    \begin{lemma}[Weak quadrature]
        For a linear flux $\mathbf{F}$, assuming the surface quadrature is accurate to degree $2k$ and that conditions \cref{eq:la_cond2,eq:la_cond3,eq:la_cond_q} are satisfied modally for some $\tb{Q}$, then let $\{\omega_i\}$ be some quadrature that is sufficiently strong. Then the condition on stability becomes
        \begin{equation}
            \mathbf{C} = \mathbf{R}_{qs}\left(\mathbf{M}_q + (\mathbf{V}^{-1}\mathbf{R}_{qk})^T\tb{Q}\mathbf{V}^{-1}\mathbf{R}_{qk}\right)^{-1}(\mathbf{L}_\partial\mathbf{P}_{qk})^T\mathbf{W}_\partial
        \end{equation}
        where $\mathbf{M}_q=\mathrm{diag}(\omega)$ and the $L_2$ projection and restriction matrices are defined as in \cref{def:l2_proj}.
    \end{lemma}
    \begin{proof}
        Starting from the statement of the FR method we have
        \begin{equation}
            \px{}{t}\mathbf{u} = -\mathbf{D}\mathbf{F} - \mathbf{C}\left[(\mathbf{n}\cdot\mathbf{F}^{num}) -\mathbf{N}\hb{L}_\partial\mathbf{F}\right].
        \end{equation}
        To integrate with sufficient accuracy we wish to use $\mathbf{M}_q + \mathbf{Q}_q$, therefore we multiply by $\mathbf{u}^T\mathbf{P}_{qk}^T(\mathbf{M}_q + \mathbf{Q}_q)\mathbf{P}_{qk}$ to obtain
        \begin{equation}
            \mathbf{u}^T\mathbf{P}_{qk}^T(\mathbf{M}_q + \mathbf{Q}_q)\mathbf{P}_{qk}\px{}{t}\mathbf{u} = -\mathbf{u}^T\mathbf{P}_{qk}^T(\mathbf{M}_q + \mathbf{Q}_q)\mathbf{P}_{qk}\mathbf{D}\mathbf{F} - \mathbf{u}^T\mathbf{P}_{qk}^T(\mathbf{M}_q + \mathbf{Q}_q)\mathbf{P}_{qk}\mathbf{C}\left[(\mathbf{n}\cdot\mathbf{F}^\mathrm{num}) -\mathbf{N}\hb{L}_\partial\mathbf{F}\right].
        \end{equation}
        As before we require
        \begin{equation}
            \mathbf{P}_{qk}^T(\mathbf{M}_q + \mathbf{Q}_q)\mathbf{P}_{qk}\mathbf{C} = \mathbf{L}_\partial^T\mathbf{W}_\partial
        \end{equation}
        for stability. By definition, we have $\mathbf{P}_{qk}^T\mathbf{M}_q\mathbf{P}_{qk}\mathbf{R}_{qk}\mathbf{M}_q^{-1}\mathbf{R}_{qk}^T=\mathbf{I}$ and hence we obtain
        \begin{equation}
            \mathbf{C} = \mathbf{R}_{qk}(\mathbf{M}_q+\mathbf{Q}_q)^{-1}\mathbf{R}_{qk}^T\mathbf{L}_\partial^T\mathbf{W}_\partial.
        \end{equation}
        Finally, the definition of $\mathbf{Q}_q=(\mathbf{V}^{-1}\mathbf{R}_{qk})^T\tb{Q}\mathbf{V}^{-1}\mathbf{R}_{qk}$ follows naturally. This concludes the proof.
    \end{proof}
    
\end{appendices}


\end{document}